\documentclass[preprint,12pt,authoryear]{elsarticle}

\usepackage{amsmath,amssymb}
\usepackage{amsthm}
\usepackage{booktabs}
\usepackage{graphicx}
\usepackage{hyperref}
\usepackage{xcolor}
\usepackage{tikz}
\usepackage{siunitx}
\usepackage[ruled,vlined,linesnumbered]{algorithm2e}
\usepackage{lineno}
\usepackage{subcaption}
\modulolinenumbers[5]

\theoremstyle{definition}

\newtheorem{remark}{Remark}

\journal{arXiv}

\begin{document}
\begin{frontmatter}

\title{The Combined Van-and-Mopeds Routing Problem\\
       with Time Windows: Formulation, Variants,\\
       and Solution Methods}

\author[ist]{Pedro~Lameiras}
\author[ist,inesc]{Alexandre~P.~Francisco\corref{cor}}
\author[dingoo]{Adriano Serrano}
\author[isel,inesc]{C\'atia~Vaz}

  \cortext[cor]{Corresponding author: \href{mailto:aplf@tecnico.ulisboa.pt}{aplf@tecnico.ulisboa.pt}}

\address[ist]{Instituto Superior T\'ecnico, Universidade de Lisboa, Lisbon, Portugal}
\address[dingoo]{Dingoo, Lisbon, Portugal}
\address[inesc]{INESC-ID Lisboa, Lisbon, Portugal}
\address[isel]{Instituto Superior de Engenharia de Lisboa, Instituto Politécnico de Lisboa, Lisbon, Portugal}

\begin{abstract}
We introduce and study the \emph{Combined Van-and-Mopeds Routing Problem
with Time Windows} (CVMRPTW), a novel variant of the heterogeneous fleet
vehicle routing problem motivated by last-mile parcel delivery in dense
European cities.
In the CVMRPTW, a single van serves as a mobile base from which on-demand
mopeds are deployed to reach customers located in areas inaccessible or
inefficient for vans.
The primary objective is to minimise the number of mopeds used while
guaranteeing that all delivery time windows are met.
We propose a Mixed-Integer Linear Programming (MILP) formulation that
captures vehicle synchronisation, capacity tracking, and time-window
compliance through 22 structured constraint families, and extend it to
three operationally motivated variants: the standard model, the
active-waiting-van variant, and the common-depot variant.
A lexicographic two-stage optimisation framework allows secondary
objectives—route duration, combined travel time, and combined
distance—to be optimised without sacrificing the primary objective.
Three graph-sparsification preprocessing techniques reduce the number of
variables and constraints by up to 25\%.
For instances beyond the practical limits of exact solvers, we develop
the Cluster-Based Combined Routing (CBCR) heuristic, which decomposes
the problem into a clustering phase, a cluster-level exact solve, and a
route-reconstruction phase.
Computational experiments on real-world instances derived from
OpenStreetMap data for Lisbon and Stuttgart show that the MILP solves
instances with up to 48 customers to feasibility within one hour, and
that the CBCR heuristic extends tractability to instances with up to 80 customers under favourable
conditions.
\end{abstract}


\end{frontmatter}


\section{Introduction}
\label{sec:intro}

Urban parcel delivery volumes in Europe have grown at unprecedented rates
over the past decade, driven by the rapid expansion of e-commerce and
changing consumer expectations for same-day or next-day delivery.
Almost three quarters of Europeans now live in cities, towns, or
suburbs~\citep{Eurostat2016}, and urban population growth shows no sign of
abating.
This concentration of demand, combined with the physical constraints of
historic European city centres---narrow streets, steep gradients, low-emission
zones, and congestion---poses fundamental challenges for last-mile logistics
operators relying exclusively on conventional delivery vans.

One operational response that has emerged in practice is the use of
\emph{heterogeneous fleets} that combine vehicles of different types and
capabilities.
In particular, the pairing of a van---which carries the full load from a
warehouse and can navigate major arterials efficiently---with smaller,
more agile vehicles such as cargo bikes, electric mopeds, or drones, has
attracted growing interest from both industry practitioners and the
academic operations research (OR) community~\citep{Crainic2016,
ctbrp}.
The van acts as a mobile depot or relay point, deploying the lighter
vehicles to serve customers in locations the van itself cannot efficiently
reach.

This paper is motivated by \textit{Dingoo}, a Lisbon-based parcel delivery
company whose operational setting exemplifies the challenge described above.
Lisbon's geography---characterised by steep hills, narrow historic lanes,
and a dense pedestrian-oriented urban fabric---makes it effectively
impossible for a delivery van to reach a significant fraction of the
customer base within the time constraints imposed by the market.
Dingoo's solution combines a single delivery van with a fleet of
on-demand electric mopeds: the van follows a coordinated route through the
city while dispatching mopeds at selected \emph{combined nodes} to serve
nearby customers.
The mopeds, once loaded at the van, operate independently until they
rejoin the van at the next rendezvous point.
The core operational objective is to fulfil all deliveries within their
time windows while deploying as few mopeds as possible, thereby minimising
the variable cost of on-demand moped riders.

Despite the richness of the vehicle routing literature, no existing
problem class fully captures this setting.
The closest relatives are the Capacitated Vehicle Routing Problem with
Time Windows~\citep[CVRPTW;][]{vrp3}, the Heterogeneous Fleet
Vehicle Routing Problem~\citep[HFVRP;][]{Gendreau1999}, and the Combined
Truck-and-Cargo-Bike Routing Problem~\citep[CTBRP;][]{ctbrp}.
The CVMRPTW inherits complexity from all three: it requires routing two
types of vehicles with different capacity and accessibility profiles,
enforcing hard time windows, and imposing synchronisation constraints
between the van and the mopeds at combined nodes.
To the best of our knowledge, the CVMRPTW has not been previously
studied in the OR literature.

\medskip
\noindent\textbf{Contributions.}
This paper makes the following contributions:

\begin{enumerate}
  \item We formally define the \emph{Combined Van-and-Mopeds Routing Problem
        with Time Windows} (CVMRPTW) and provide its first Mixed-Integer
        Linear Programming formulation, comprising 22 structured constraint
        families covering flow conservation, synchronisation, capacity
        tracking, and time-window compliance.

  \item We introduce three operationally motivated problem variants---the
        standard model~(\textsc{s}), the active-waiting-van
        model~(\textsc{awv}), and the common-depot model~(\textsc{cd})---and
        a lexicographic two-stage optimisation framework that supports four
        selectable secondary objectives.

  \item We propose three graph-sparsification preprocessing techniques that
        reduce the number of decision variables and constraints by up to
        25\%, improving solver performance without loss of optimality.

  \item We develop the \emph{Cluster-Based Combined Routing} (CBCR)
        heuristic, a decomposition approach that makes the problem tractable
        for instances beyond the practical limits of exact solvers, and that
        is adaptable to related subclasses of the VRP.

  \item We construct a real-world instance generator using OpenStreetMap
        data for Lisbon and Stuttgart, introduce new performance metrics
        tailored to the CVMRPTW (Average Moped Load, Average Synchronisation
        Ratio), and conduct extensive computational experiments. Code, tests and results
        are available at \url{https://doi.org/10.5281/zenodo.20667375}.
\end{enumerate}

\medskip
The remainder of the paper is organised as follows.
Section~\ref{sec:literature} reviews the relevant literature and
positions the CVMRPTW relative to existing problem classes.
Section~\ref{sec:problem} formally defines the problem, introduces the
graph model, and presents the MILP formulation.
Section~\ref{sec:variants} describes the three problem variants and the
lexicographic optimisation framework.
Section~\ref{sec:preprocessing} presents the graph-sparsification
preprocessing techniques.
Section~\ref{sec:heuristic} introduces the CBCR heuristic.
Section~\ref{sec:experiments} reports computational experiments; we relied on 
Gurobi 12.0.2~\citep{gurobi} for solving MILP instances.
Section~\ref{sec:conclusions} concludes with a summary of findings,
managerial insights, and directions for future research.

\section{Literature review}
\label{sec:literature}

The CVMRPTW sits at the intersection of several well-studied combinatorial
optimisation problems.
This section reviews the relevant strands of the vehicle routing
literature and closes by positioning the CVMRPTW with respect to existing
problem classes.

\subsection{Vehicle routing fundamentals}
\label{subsec:lit:vrp}

The \emph{Travelling Salesman Problem} (TSP) is the canonical combinatorial
optimisation problem: given a set of cities and pairwise travel costs, find
the minimum-cost Hamiltonian cycle.
Its NP-hardness was established by \citet{Karp1972}, and it remains a
benchmark for exact, heuristic, and metaheuristic methods alike.
The \emph{Miller--Tucker--Zemlin} (MTZ) formulation \citep{Miller1960}
provides a compact MILP encoding with a polynomial number of constraints,
at the cost of a weaker linear relaxation compared to subtour-elimination
cuts; this trade-off is central to the model design in the present paper.

The \emph{Capacitated Vehicle Routing Problem} (CVRP) generalises the TSP
to a fleet of $K$ identical vehicles of capacity $C$, each originating from
a central depot, that must collectively serve all customers while minimising
total routing cost \citep{vrp3}.
The CVRP is NP-hard by reduction from TSP \citep{Karp1972}, and its
complexity grows rapidly with the number of customers and vehicles.

The \emph{Capacitated VRP with Time Windows} (CVRPTW) extends the CVRP by
associating each customer $i$ with a hard time window $[a_i, b_i]$: a
vehicle arriving before $a_i$ must wait, and service must begin no later
than $b_i$.
The CVRPTW is also NP-hard by restriction \citep{vrp3}, and the
interaction between capacity and time constraints makes it substantially
harder in practice.
The seminal exact branch-and-bound approach of \citet{Solomon1987} and
the benchmark instances he introduced remain reference points for the
field.
State-of-the-art exact methods rely on branch-and-price, exploiting the
structure of column generation over routes \citep{Desrochers1992}, while
metaheuristics such as adaptive large neighbourhood search
\citep{Ropke2006} achieve near-optimal solutions on instances with
hundreds of customers.

\subsection{Heterogeneous fleet vehicle routing}
\label{subsec:lit:hfvrp}

In many real-world distribution settings, the vehicle fleet is not
homogeneous: different vehicle types have different load capacities,
operating costs, and accessibility profiles.
The \emph{Heterogeneous Fleet VRP} (HFVRP) relaxes the assumption of
identical vehicles, partitioning the fleet into vehicle types $p \in P$,
each characterised by capacity $C_p$, variable cost $c^p_{ij}$ per
traversed arc, and fixed activation cost $\mathit{FC}_p$
\citep{vrp3,Gendreau1999}.
Two main objective variants are distinguished: minimising total routing
cost (fixed plus variable) and minimising the number of vehicles used.
The latter is directly relevant to the CVMRPTW, where the primary
objective is to minimise the number of mopeds deployed.

\citet{Gendreau1999} proposed an early and influential tabu search
heuristic for the HFVRP, demonstrating that mixed fleets could be handled
effectively without exact enumeration.
Subsequent work has extended the HFVRP to incorporate time windows
\citep{Liu1999}, multiple depots \citep{Salhi2014}, and environmental
objectives \citep{Kocc2016}.
The HFVRP with time windows is NP-hard and constitutes one of the three
complexity sources inherited by the CVMRPTW.

\subsection{Two-echelon and collaborative vehicle routing}
\label{subsec:lit:two_echelon}

A distinct but related research direction concerns \emph{two-echelon}
or \emph{collaborative} routing, in which a primary vehicle operates as
a mobile relay for one or more secondary vehicles.
This paradigm has gained considerable traction in urban last-mile
logistics, where access restrictions and congestion prevent a single
vehicle type from efficiently serving the entire customer set
\citep{Crainic2016}.

The \emph{Two-Echelon VRP} (2E-VRP) uses fixed satellite depots as
intermediate transfer points between large freight vehicles and smaller
city vehicles \citep{Perboli2011}.
Satellites are fixed locations, known a priori, and do not move during
the planning horizon.
\citet{anderluh} extend the 2E-VRP to a multi-objective setting with
vehicle synchronisation and a ``grey zone'' where both vehicle types can
deliver, solved via large neighbourhood search; this is the closest
antecedent to the CVMRPTW in terms of synchronisation structure,
although combined nodes are fixed satellites rather than dynamically
chosen stops on the van's route.

More recently, the idea of a \emph{mobile} relay has been formalised in
the \emph{Truck-and-Drone} routing literature
\citep{MurrayC2015,Agatz2018}, where an unmanned aerial vehicle is
launched and recovered from a moving truck.
These models share the synchronisation requirement of the CVMRPTW
(secondary vehicle must meet the primary vehicle at a rendezvous point)
but differ in vehicle capabilities, capacity structure, and the absence
of hard time windows at customers in many formulations.

The \emph{Combined Truck-and-Cargo-Bike Routing Problem}
(CTBRP)~\citep{ctbrp} is the closest existing problem to the CVMRPTW.
In the CTBRP, a single truck and a single cargo bike share a set of
customers; some customers (\emph{combined nodes}) are served by the
truck and used as reloading points for the bike.
\citet{ctbrp} provide a MILP formulation and solve instances with up to
25 customers exactly.
The CVMRPTW generalises the CTBRP in three important directions: (i) the
secondary fleet consists of \emph{multiple} mopeds rather than a single
bike; (ii) \emph{hard time windows} are imposed on all customers; and
(iii) the number of secondary vehicles is itself a decision variable
subject to minimisation.

\subsection{Inter-vehicle synchronisation in vehicle routing}
\label{subsec:lit:sync}

Synchronisation between vehicles---requiring that certain events at
shared locations happen simultaneously or in a prescribed
order---introduces a qualitatively new layer of difficulty in routing
models \citep{Drexl2012}.
\citet{Drexl2012} surveys inter-vehicle synchronisation constraints
across several VRP variants, distinguishing \emph{operation}
synchronisation (vehicles must be at the same location at the same time),
\emph{load} synchronisation (transfer of goods between vehicles), and
\emph{movement} synchronisation (vehicles must depart together).
The CVMRPTW involves all three: mopeds must arrive at a combined node
before the van (operation), reload their cargo from the van (load), and
then operate independently until the next combined node (movement).

Time-window propagation provides an effective implicit mechanism for
synchronisation in MILP formulations: by encoding arrival times as
continuous variables and linking them through big-$M$ constraints
activated by binary routing variables, subtour elimination and
synchronisation can both be achieved without the MTZ auxiliary variables
\citep{vrp3,ctbrp}.
This approach is adopted in the CVMRPTW formulation developed in
Section~\ref{sec:problem}.

\subsection{Positioning the CVMRPTW}
\label{subsec:lit:positioning}

Table~\ref{tab:lit} summarises the key structural features of the
problem classes reviewed above and positions the CVMRPTW among them.
The CVMRPTW is the first problem that simultaneously combines a
\emph{variable-size secondary fleet}, \emph{hard time windows}, \emph{van-inaccessible customer zones}, and \emph{dynamic combined-node
selection} within a single MILP formulation.
No existing benchmark dataset or exact-solution methodology covers this
combination; the contributions of the present paper therefore fill a
genuine gap in the vehicle routing literature.

\begin{table}[!tbp]
  \centering
  \caption{Comparison of related vehicle routing problem classes.
           \checkmark~= feature present; ---~= feature absent;
           ($\star$)~= present in some variants only.}
  \label{tab:lit}
  \small
  \begin{tabular}{lcccccc}
    \toprule
    & \rotatebox{66}{\textbf{CVRPTW}} & \rotatebox{66}{\textbf{HFVRP}} & \rotatebox{66}{\textbf{2E-VRP}}
    & \rotatebox{66}{\textbf{Truck-Drone}} & \rotatebox{66}{\textbf{CTBRP}} & \rotatebox{66}{\textbf{CVMRPTW}} \\
    \midrule
    Heterogeneous fleet       & ---        & \checkmark & \checkmark & \checkmark & \checkmark & \checkmark \\
    Hard time windows         & \checkmark & ($\star$)  & ($\star$)  & ---        & ---        & \checkmark \\
    Mobile relay / rendezvous & ---        & ---        & ---        & \checkmark & \checkmark & \checkmark \\
    Inter-vehicle sync        & ---        & ---        & \checkmark & \checkmark & \checkmark & \checkmark \\
    Variable fleet size (min) & ---        & ($\star$)  & ---        & ---        & ---        & \checkmark \\
    Dynamic combined nodes    & ---        & ---        & ---        & \checkmark & ---        & \checkmark \\
    Multiple secondary veh.   & ---        & \checkmark & \checkmark & ($\star$)  & ($\star$)  & \checkmark \\
    \bottomrule
  \end{tabular}
\end{table}

\section{Problem definition and MILP formulation}
\label{sec:problem}

\subsection{Motivating setting and informal description}
\label{subsec:prob:motivation}

The CVMRPTW is motivated by the real-world operations of \textit{Dingoo},
a parcel delivery company serving the city of Lisbon, Portugal.
Lisbon's geography presents exceptional last-mile challenges: steep
hillside neighbourhoods such as Alfama and Bairro Alto are threaded by
cobblestone alleys and stairways that are inaccessible to delivery vans,
while pedestrian zones and strict low-emission regulations further
restrict vehicular access throughout the historic centre.
At the same time, the city hosts a dense network of food-and-beverage
businesses---cafés, restaurants, delicatessens---that require frequent
restocking with time-sensitive perishable goods.

Dingoo's operational model is designed to cope with this environment.
A single delivery van loads packages at a central warehouse
(located in Azambuja or Sacavém, approximately 45~km from the city
centre) and follows an optimised route through Lisbon.
At selected \emph{combined nodes} along its route, the van serves as a
mobile loading station: electric mopeds, deployed on demand, collect
packages from the van and complete the last segment of delivery to nearby
customers that the van cannot reach.
Mopeds operate independently after each reload and may terminate their
route at any customer location; they do not need to return to a fixed
depot.
All deliveries must respect hard time windows, and the number of mopeds
activated per shift---charged at an hourly on-demand rate---must be
minimised.

Operations are organised into three daily shifts.
The van carries up to approximately 50 packages per shift; each moped
carries up to $L_m$ packages (typically 4--6).
Time windows are contractually binding and define the service-level
agreement with customers.
The van's route is constrained by the LIFO (last-in, first-out) loading
order of its package grid, leaving little flexibility for reordering
deliveries after departure.

\subsection{Graph model and notation}
\label{subsec:prob:graph}

The distribution network is represented as a directed weighted graph
$G = (V, E)$.
The node set $V$ contains a depot $v_s$ (van start), an end node $v_e$
(van terminus), and $n$ customer nodes $v_1, v_2, \dots, v_n$.
We write $V_c = V \setminus \{v_s, v_e\}$ for the set of customer nodes.

Not all nodes are accessible to every vehicle type.
Let $V_{vn} \subseteq V_c$ be the set of customers reachable by the van,
and $V_m \subseteq V_c$ the set reachable by mopeds.
A customer $v \in V_c \setminus V_{vn}$ lies in a van-restricted zone
and must be served exclusively by a moped.
Conversely, a customer $v \in V_c \setminus V_m$ (rare in practice) must
be served by the van.

The edge set $E$ is split into a van subgraph $E_{vn}$ and a moped
subgraph $E_m$ after a preprocessing step (Section~\ref{sec:preprocessing}).
For each edge $(i,j) \in E_{vn}$ the van travel time is $c^{vn}_{(i,j)}$,
and the van-distance is $c'^{vn}_{(i,j)}$.
Analogously, $c^m_{(i,j)}$ and $c'^m_{(i,j)}$ are the moped travel time
and distance for $(i,j) \in E_m$.
Each customer $v \in V_c$ has a positive demand $d_v > 0$; the depot and
end node have zero demand.
Each customer $v$ also has a hard time window $[a_v, b_v]$ and vehicle-type
specific service times $s^{vn}_v$ and $s^m_v$ (time required to
complete the delivery once the vehicle arrives).

A node $v \in V_c$ is a \emph{combined node} if it is visited by both
the van and a moped: the van delivers its own package(s) there and
simultaneously reloads the moped.
The set of combined nodes $K \subseteq V_c$ is not fixed a priori; it
emerges as part of the optimisation.

Table~\ref{tab:notation} summarises the main sets, parameters, and
decision variables used throughout the formulation.

\begin{table}[!tbp]
  \centering
  \caption{Summary of notation.}
  \label{tab:notation}
  \small
  \begin{tabular}{ll}
    \toprule
    \textbf{Symbol} & \textbf{Meaning} \\
    \midrule
    \multicolumn{2}{l}{\textit{Sets}} \\
    $V$                   & All nodes: $\{v_s\} \cup V_c \cup \{v_e\}$ \\
    $V_c$                 & Customer nodes \\
    $V_{vn},\; V_m$       & Van-accessible / moped-accessible customers \\
    $E_{vn},\; E_m$       & Van / moped edge sets \\
    $K$                   & Combined nodes (endogenous) \\
    \midrule
    \multicolumn{2}{l}{\textit{Parameters}} \\
    $d_v$                 & Demand at node $v$ \\
    $[a_v, b_v]$          & Hard time window at node $v$ \\
    $s^{vn}_v,\; s^m_v$   & Service time for van / moped at node $v$ \\
    $c^{vn}_{(i,j)},\; c^m_{(i,j)}$ & Van / moped travel time on edge $(i,j)$ \\
    $c'^{vn}_{(i,j)},\; c'^m_{(i,j)}$ & Van / moped travel distance on edge $(i,j)$ \\
    $L_m$                 & Moped capacity \\
    $M$                   & Big-$M$ constant: $M = T_{\mathrm{shift}}$ (shift duration)\\
    \midrule
    \multicolumn{2}{l}{\textit{Decision variables}} \\
    $x^{vn}_v \in \{0,1\}$     & 1 if van visits node $v$ \\
    $x^m_v \in \{0,1\}$        & 1 if a moped visits node $v$ \\
    $x^{vn}_{(i,j)} \in \{0,1\}$ & 1 if van traverses edge $(i,j)$ \\
    $x^m_{(i,j)} \in \{0,1\}$  & 1 if a moped traverses edge $(i,j)$ \\
    $T^{vn}_v \geq 0$          & Van arrival time at node $v$ \\
    $T^m_v \geq 0$             & Moped arrival time at node $v$ \\
    $y^m \in \mathbb{N}_0$     & Number of mopeds used (primary objective) \\
    $z^m_v \geq 0$             & Auxiliary: net moped outflow at node $v$ \\
    $\ell_v \geq 0$            & Cumulative moped load since last reload \\
    \bottomrule
  \end{tabular}
\end{table}

\subsection{MILP formulation}
\label{subsec:prob:milp}

The primary objective $f_1$ is to minimise the number of mopeds deployed:
\begin{equation}
  \min \quad y^m. \label{obj:primary}
\end{equation}
The 22 constraint families are presented below, grouped by function.

\medskip
\noindent\textbf{C1 -- Van origin.}
Only the van departs from the depot; mopeds do not:
\begin{equation}
  x^{vn}_{v_s} = 1, \quad x^m_{v_s} = 0. \label{c1}
\end{equation}

\medskip
\noindent\textbf{C2 -- Van terminus.}
Only the van terminates at the designated end node, which cannot be a
combined node:
\begin{equation}
  x^{vn}_{v_e} = 1,\quad x^m_{v_e} = 0,\quad
  \sum_{(v_e,\,w)\in E_{vn}} x^{vn}_{(v_e,w)} = 0. \label{c2}
\end{equation}

\medskip
\noindent\textbf{C3 -- Time origin.}
The clock starts when the van departs from the depot:
\begin{equation}
  T^{vn}_{v_s} = 0. \label{c3}
\end{equation}

\medskip
\noindent\textbf{C4--C5 -- Vehicle-accessibility.}
The van is forbidden from visiting nodes outside $V_{vn}$; mopeds are
forbidden from visiting nodes outside $V_m$:
\begin{align}
  x^{vn}_v &= 0, \quad \forall v \in V_c \setminus V_{vn}, \label{c4} \\
  x^m_v    &= 0, \quad \forall v \in V_c \setminus V_m.   \label{c5}
\end{align}
After graph sparsification (Section~\ref{sec:preprocessing}) these
constraints are satisfied implicitly and can be omitted.

\medskip
\noindent\textbf{C6 -- Coverage.}
Each customer must be served by at least one vehicle:
\begin{equation}
  x^{vn}_v + x^m_v \geq 1, \quad \forall v \in V_c. \label{c6}
\end{equation}

\medskip
\noindent\textbf{C7--C8 -- Van flow conservation.}
If the van visits a node, exactly one incoming and one outgoing van
edge must be active:
\begin{align}
  x^{vn}_v &\leq \sum_{(w,v)\in E_{vn}} x^{vn}_{(w,v)},
    \quad \forall v \in V,\; v \neq v_s, \label{c7} \\
  \sum_{(v,w)\in E_{vn}} x^{vn}_{(v,w)} &= x^{vn}_v,
    \quad \forall v \in V,\; v \neq v_e. \label{c8}
\end{align}

\medskip
\noindent\textbf{C9 -- Moped count.}
The number of mopeds used equals the total net outflow of mopeds from
combined nodes.
To avoid the bilinear product $x^m_{(v,w)} \cdot x^{vn}_v$, a
continuous auxiliary variable $z^m_v \geq 0$ is introduced and the
constraint is linearised:
\begin{align}
  \sum_{(v,w)\in E_m} x^m_{(v,w)} - \sum_{(w,v)\in E_m} x^m_{(w,v)}
    &\leq z^m_v + (1 - x^{vn}_v)\cdot M,
    \quad \forall v \in V_c, \label{c9a} \\
  \sum_{v \in V_c} z^m_v &\leq y^m. \label{c9b}
\end{align}

\medskip
\noindent\textbf{C10 -- Moped flow at combined nodes.}
Mopeds always depart from combined nodes; they cannot simply pass
through, and the moped flow between two consecutive combined nodes
must be zero:
\begin{align}
  \sum_{(w,v)\in E_m} x^m_{(w,v)} &\leq
    \sum_{(v,w)\in E_m} x^m_{(v,w)} + (1 - x^{vn}_v)\cdot M,
    \quad \forall v \in V_c, \label{c10a} \\
  x^{vn}_v + x^{vn}_w &\leq 1 + (1 - x^m_{(v,w)})\cdot M,
    \quad \forall (v,w) \in E_m. \label{c10b}
\end{align}

\medskip
\noindent\textbf{C11 -- Moped flow at non-combined nodes.}
At moped-only nodes exactly one incoming moped edge must be active:
\begin{align}
  1 &\leq \sum_{(w,v)\in E_m} x^m_{(w,v)} + x^{vn}_v \cdot M,
    \quad \forall v \in V_c, \label{c11a} \\
  \sum_{(w,v)\in E_m} x^m_{(w,v)} &\leq 1 + x^{vn}_v \cdot M,
    \quad \forall v \in V_c. \label{c11b}
\end{align}

\medskip
\noindent\textbf{C12 -- Moped start restriction.}
Mopeds cannot start their shift at non-combined nodes; outgoing moped
edges at non-combined nodes are bounded to at most one:
\begin{equation}
  \sum_{(v,w)\in E_m} x^m_{(v,w)} \leq 1 + x^{vn}_v \cdot M,
  \quad \forall v \in V_c. \label{c12}
\end{equation}

\medskip
\noindent\textbf{C13 -- Moped edge–node consistency.}
Active moped edges may only connect visited nodes:
\begin{align}
  x^m_{(w,v)} &\leq x^m_v, \quad \forall (w,v) \in E_m, \label{c13a} \\
  x^m_{(w,v)} &\leq x^m_w, \quad \forall (w,v) \in E_m. \label{c13b}
\end{align}

\medskip
\noindent\textbf{C14 -- Moped load initialisation.}
The load tracker $\ell_v$ is zero at combined nodes (van-delivered) and
bounded by $L_m$ at moped-only nodes:
\begin{equation}
  \ell_v \leq (1 - x^{vn}_v) \cdot L_m, \quad \forall v \in V. \label{c14}
\end{equation}

\medskip
\noindent\textbf{C15 -- Moped capacity tracking.}
The cumulative moped load increases along each moped sub-route; upon
reaching a combined node the counter resets:
\begin{equation}
  \ell_v + d_w \leq \ell_w +
    \bigl(1 - x^m_{(v,w)} + x^{vn}_w\bigr)
    \cdot \bigl(L_m + \max_{v' \in V} d_{v'}\bigr),
  \quad \forall (v,w) \in E_m. \label{c15}
\end{equation}

\medskip
\noindent\textbf{C16 -- Time-variable unification at non-combined nodes.}
For nodes visited by only one vehicle type the two time variables
coincide (linearised via big-$M$):
\begin{align}
  T^{vn}_v &\leq T^m_v + (x^{vn}_v + x^m_v - 1)\cdot M,
    \quad \forall v \in V, \label{c16a} \\
  T^m_v    &\leq T^{vn}_v + (x^{vn}_v + x^m_v - 1)\cdot M,
    \quad \forall v \in V. \label{c16b}
\end{align}

\medskip
\noindent\textbf{C17 -- Time propagation and synchronisation.}
Arrival times must be consistent with travel times, service times, and
the synchronisation requirement at combined nodes.
For the van sub-graph:
\begin{equation}
  T^{vn}_v + s^{vn}_v + c^{vn}_{(v,w)}
  \leq T^{vn}_w + (1 - x^{vn}_{(v,w)})\cdot M,
  \quad \forall (v,w) \in E_{vn}. \label{c17van}
\end{equation}
For the moped sub-graph the synchronisation factor
$(T^{vn}_v - T^m_v)$ is added at combined nodes (it equals zero
elsewhere by C16):
\begin{multline}
  T^m_v + (T^{vn}_v - T^m_v)
  + (1 - x^{vn}_v)\cdot s^m_v + x^{vn}_v \cdot s^{vn}_v
  + c^m_{(v,w)} \\
  \leq T^m_w + (1 - x^m_{(v,w)})\cdot M,
  \quad \forall (v,w) \in E_m. \label{c17mop}
\end{multline}
Time propagation from C17 implicitly eliminates subtours; no
Miller--Tucker--Zemlin constraints are required.

\medskip
\noindent\textbf{C18 -- Van time windows.}
Arrival times of the van must respect (combined) nodes time windows:
\begin{equation}
  a_v \leq T^{vn}_v \leq b_v,
  \quad \forall v \in V_{vn}. \label{c18}
\end{equation}

\medskip
\noindent\textbf{C19 -- Moped time windows.}
Time windows apply also to moped-only nodes (linearised via big-$M$):
\begin{align}
  a_v &\leq T^m_v + (x^{vn}_v + x^m_v - 1)\cdot M,
    \quad \forall v \in V_m, \label{c19a} \\
  T^m_v &\leq b_v + (x^{vn}_v + x^m_v - 1)\cdot M,
    \quad \forall v \in V_m. \label{c19b}
\end{align}

\medskip
\noindent\textbf{C20 -- Moped-before-van synchronisation.}
At every combined node the moped must arrive before the van:
\begin{equation}
  T^m_v \leq T^{vn}_v + (2 - x^{vn}_v - x^m_v)\cdot M,
  \quad \forall v \in V. \label{c20}
\end{equation}

\medskip
\noindent\textbf{C21 -- Binary variables.}
\begin{align}
  x^{vn}_v,\; x^m_v &\in \{0,1\},\; \forall v \in V,\\
  x^{vn}_{(i,j)} &\in \{0,1\},\; \forall (i,j) \in E_{vn},\\
  x^m_{(i,j)} &\in \{0,1\},\; \forall (i,j) \in E_m. \label{c21}
\end{align}

\medskip
\noindent\textbf{C22 -- Continuous and integer variables.}
\begin{align}
  T^{vn}_v,\; T^m_v &\in \mathbb{R}_{\geq 0},\; \forall v \in V,\\
  y^m,\; z^m_v,\; \ell_v &\in \mathbb{N}_0,\; \forall v \in V. \label{c22}
\end{align}

\subsection{Big-\textit{M} value}
\label{subsec:prob:bigM}

The constant $M$ must dominate all time and flow quantities in the
model.
Since the shift duration $T_{\mathrm{shift}}$ upper-bounds all arrival
times, and the number of customers served per shift is at most
$|V_c| \leq 50$ in the motivating setting (bounded by van capacity),
we set
\begin{equation}
  M = \max\bigl(|V| - 1,\; T_{\mathrm{shift}}\bigr) = T_{\mathrm{shift}},
\end{equation}
because $T_{\mathrm{shift}} \gg |V|$ in all practical instances.
This tight choice improves the LP relaxation bound and thereby the
branch-and-bound performance.

\subsection{Model complexity}
\label{subsec:prob:complexity}

The number of binary variables scales as $O(n^2)$ (from the edge sets
$E_{vn}$ and $E_m$) and the number of continuous variables as $O(n)$.
The number of constraints is dominated by C15 and C17, both indexed
over edges, giving $O(n^2)$ constraints in total.
Since the CVRPTW is NP-hard~\citep{vrp3} and reduces to the CVMRPTW
by setting $V_m = \emptyset$ (forcing all deliveries via the van), the
CVMRPTW is also NP-hard.
The combination of capacity (C14--C15), time-window (C18--C19), and
synchronisation (C20) constraints, together with the endogenous fleet
size objective, makes it strictly harder in practice than any single
one of its constituent problem classes.

\section{Problem variants and lexicographic optimisation}
\label{sec:variants}

The ground model of Section~\ref{sec:problem} captures Dingoo's
baseline operations but admits several natural extensions that reflect
alternative operational policies.
Three variants are introduced below.
They share the same graph structure and the majority of the constraint
set; only the depot-assignment and synchronisation rules differ.
A lexicographic two-stage optimisation framework applicable to all
variants closes the section.

\subsection{Standard variant (\textsc{s})}
\label{subsec:var:standard}

The standard variant corresponds to the ground model of
Section~\ref{sec:problem} exactly.
Mopeds do not share a depot with the van: they are stationed throughout
the city and are dispatched on demand.
Consequently, every moped route must begin at a combined node, i.e.\ at
a point where the van is present to load it.
Constraint~\eqref{c1} enforces this by forbidding $x^m_{v_s} = 1$.
The synchronisation requirement of constraint~\eqref{c20} guarantees that each
moped arrives at its departure combined node before the van does,
ensuring cargo is available when the moped is activated.

This variant best reflects Dingoo's current operational setup and is
used as the reference configuration in all computational comparisons.

\subsection{Active-waiting-van variant (\textsc{awv})}
\label{subsec:var:awv}

In the standard variant the van delivers at a combined node and reloads
the moped in a single synchronised event.
The \textsc{awv} variant relaxes this by decoupling the two events: a
moped may be reloaded by the van \emph{before} the van's own scheduled
delivery at that node.
If the van arrives before the opening of the time window $a_v$, it
waits and serves the customer at $a_v$; the moped, however, can be
loaded and dispatched as soon as the van arrives, gaining a potential
head-start.

This is achieved by three targeted changes to the ground model.
First, constraints~\eqref{c17van}--\eqref{c17mop} are replaced to allow the van's
effective departure to be $\max(T^{vn}_v, a_v)$ at combined nodes:
\begin{align}
  \max(T^{vn}_v,\,a_v) + s^{vn}_v + c^{vn}_{(v,w)}
    &\leq T^{vn}_w + (1 - x^{vn}_{(v,w)})\cdot M,
    \quad \forall (v,w)\in E_{vn}, \label{awv:van} \\
  T^m_v + \bigl(\max(T^{vn}_v,\,T^m_v) & - T^m_v\bigr)
    + (1-x^{vn}_v)\cdot s^m_v \nonumber \\
    +\ x^{vn}_v\cdot s^{vn}_v + c^m_{(v,w)}
    &\leq T^m_w + (1 - x^m_{(v,w)})\cdot M,
    \quad \forall (v,w)\in E_m. \label{awv:mop}
\end{align}
The synchronisation factor in inequality \eqref{awv:mop} is now
$\max(T^{vn}_v, T^m_v) - T^m_v \geq 0$, so the moped waits only if the
van has not yet arrived---the reverse never occurs by construction.
Second, the van time-window lower bound in constraint~\eqref{c18} is relaxed to
allow the recorded arrival time to precede $a_v$:
\begin{equation}
  T^{vn}_v \leq b_v, \quad \forall v \in V_{vn}. \label{awv:tw}
\end{equation}
Third, constraint~\eqref{c20} is replaced to allow the moped to arrive before
$\max(T^{vn}_v, a_v)$:
\begin{equation}
  T^m_v \leq \max(T^{vn}_v,\,a_v) + (2-x^{vn}_v-x^m_v)\cdot M,
  \quad \forall v \in V. \label{awv:sync}
\end{equation}
The $\max$ operator is handled natively by Gurobi via
\texttt{addGenConstrMax}, avoiding auxiliary variables.
Note that formulating the synchronisation factor as
$\max(T^{vn}_v - T^m_v, 0)$ would introduce floating-point
subtraction instability; the form in inequality \eqref{awv:mop} is numerically
safer.

\subsection{Common-depot variant (\textsc{cd})}
\label{subsec:var:cd}

In the \textsc{cd} variant all mopeds depart from the same depot as the
van ($v_s$), reflecting a scenario in which all vehicles are
pre-loaded at the warehouse before the shift begins.
This eliminates the combined-node-as-start-point requirement and allows
mopeds to fan out across the city independently from the first moment.

The structural changes to the ground model are as follows.
Constraint~\eqref{c1} is replaced by:
\begin{equation}
  x^{vn}_{v_s} = 1, \quad x^m_{v_s} = 1, \label{cd:c1}
\end{equation}
and constraint~\eqref{c3} by:
\begin{equation}
  T^{vn}_{v_s} = 0, \quad T^m_{v_s} = 0. \label{cd:c3}
\end{equation}
The moped flow constraints \eqref{c7}--\eqref{c13a} (excepting the
second sub-constraint of C10) are replaced by four dedicated families
that reflect the common-depot topology.

\medskip
\noindent\textbf{C25 -- Van flow (equality).}
Incoming and outgoing van edges at each visited node are exactly one:
\begin{equation}
  \sum_{(w,v)\in E_{vn}} x^{vn}_{(w,v)} = x^{vn}_v,
  \quad \forall v \in V,\; v \neq v_s. \label{cd:c25}
\end{equation}

\medskip
\noindent\textbf{C26 -- Van outflow.}
\begin{equation}
  \sum_{(v,w)\in E_{vn}} x^{vn}_{(v,w)} = x^{vn}_v,
  \quad \forall v \in V,\; v \neq v_e. \label{cd:c26}
\end{equation}

\medskip
\noindent\textbf{C27 -- Moped flow balance.}
Mopeds may end their route at any customer node; incoming flow must
weakly exceed outgoing flow at all non-depot nodes:
\begin{align}
  \sum_{(w,v)\in E_m} x^m_{(w,v)}
    &\geq \sum_{(v,w)\in E_m} x^m_{(v,w)},
    \quad \forall v \in V,\; v \neq v_s, \label{cd:c27a} \\
  \sum_{(w,v)\in E_m} x^m_{(w,v)}
    &\leq \sum_{(v,w)\in E_m} x^m_{(v,w)} + (2-x^{vn}_v-x^m_v)\cdot M,
    \quad \forall v \in V,\; v \neq v_s. \label{cd:c27b}
\end{align}

\medskip
\noindent\textbf{C28 -- Moped visit bounds.}
At each visited node the number of incoming moped edges is at least one
and at most $y^m$; no two mopeds may terminate at the same node
(required for the \texttt{tct} secondary objective to be well-defined):
\begin{align}
  x^m_v &\leq \sum_{(w,v)\in E_m} x^m_{(w,v)} \leq x^m_v \cdot y^m,
    \quad \forall v \in V,\; v \neq v_s, \label{cd:c28a} \\
  \sum_{(w,v)\in E_m} x^m_{(w,v)} &\leq 1 + x^{vn}_v \cdot M,
    \quad \forall v \in V,\; v \neq v_s. \label{cd:c28b}
\end{align}

\medskip
\noindent\textbf{C29 -- Moped departure count.}
Exactly $y^m$ mopeds depart from the depot:
\begin{equation}
  \sum_{(v_s,\,w)\in E_m} x^m_{(v_s,w)} = y^m. \label{cd:c29}
\end{equation}
The auxiliary variable $z^m_v$ is no longer needed in the \textsc{cd}
variant and is removed from constraint~\eqref{c22}.

\subsection{Lexicographic two-stage optimisation}
\label{subsec:var:lexico}

Minimising the number of mopeds $y^m$ is the primary objective $f_1$ for all
variants.
Once $y^m$ is determined, a secondary objective $f_2$ can be optimised
without sacrificing the primary goal.
This is formalised as a lexicographic two-stage problem:
\begin{align}
  \text{Stage 1:}\quad & y^m_* = \arg\min \; y^m
    \quad \text{subject to the CVMRPTW constraints,} \label{lex:s1} \\
  \text{Stage 2:}\quad & \min \; f_2
    \quad \text{subject to the CVMRPTW constraints and}\;
    y^m = y^m_*. \label{lex:s2}
\end{align}
Four secondary objectives $f_2$ are supported, defined in
Table~\ref{tab:objectives}.
They exploit the time and distance variables already present in the
model, so no new structural constraints are needed beyond fixing $y^m$.

\begin{table}[!tbp]
  \centering
  \caption{Secondary objective functions $f_2$. Here $\mathcal{V}$
           denotes the van route and $\mathcal{M}$ the set of moped
           routes; $T^{vn}_{v_e}$ is the van's arrival at the end node;
           $T_f(\mathcal{M})$ and $T_i(\mathcal{M})$ are the sum of
           moped route end and start times, respectively.}
  \label{tab:objectives}
  \small
  \begin{tabular}{llp{7cm}}
    \toprule
    \textbf{Code} & \textbf{Name} & \textbf{Expression} \\
    \midrule
    \texttt{vrd} & Van route duration &
      $T(\mathcal{V}) = T^{vn}_{v_e}$ \\[4pt]
    \texttt{tct} & Total completion time &
      $T(\mathcal{V}) + T_f(\mathcal{M}) - T_i(\mathcal{M})$ \\[4pt]
    \texttt{cdu} & Combined route duration &
      $\displaystyle\sum_{(v,w)\in E_{vn}} x^{vn}_{(v,w)}\,c^{vn}_{(v,w)}
       + \sum_{(v,w)\in E_m} x^m_{(v,w)}\,c^m_{(v,w)}$ \\[6pt]
    \texttt{cdi} & Combined route distance &
      $\displaystyle\sum_{(v,w)\in E_{vn}} x^{vn}_{(v,w)}\,c'^{vn}_{(v,w)}
       + \sum_{(v,w)\in E_m} x^m_{(v,w)}\,c'^m_{(v,w)}$ \\
    \bottomrule
  \end{tabular}
\end{table}

Computing $T_f(\mathcal{M})$ requires identifying the last node of each
moped route.
A binary indicator $ls_v \in \{0,1\}$ equals 1 if $v$ is the terminal
node of some moped route (i.e.\ $x^m_v = 1$ and no outgoing moped edge
is active).
An auxiliary variable $mt_v \geq 0$ captures the arrival time at $v$
when $ls_v = 1$ and is zero otherwise:
\begin{equation}
  T_f(\mathcal{M}) = \sum_{v \in V} mt_v, \quad
  mt_v \leq T^m_v, \quad
  mt_v \geq T^m_v - (1-ls_v)\cdot M, \quad
  mt_v \leq ls_v \cdot M, \label{tct:end}
\end{equation}
with $ls_v \leq x^m_v$,
$\sum_{(v,w) \in E_m} x^m_{(v,w)} \leq (1 - ls_v) \cdot M$,
and $ls_v \geq x^m_v - \sum_{(v,w)\in E_m} x^m_{(v,w)}$.
The moped start-time aggregate $T_i(\mathcal{M})$ is computed via:
\begin{equation}
  T_i(\mathcal{M}) = \sum_{v \in V} ms_v \cdot T^m_v, \label{tct:start}
\end{equation}
where $ms_v = x^{vn}_v \cdot \bigl(\sum_{(v,w)\in E_m} x^m_{(v,w)} -
\sum_{(w,v)\in E_m} x^m_{(w,v)}\bigr) \in \mathbb{N}_0$ counts the
net number of mopeds starting at combined node $v$.

\begin{remark}
  The \texttt{tct} objective jointly penalises late van completion and
  long moped routes, incentivising the solver to keep both vehicles
  active as much as possible.
  In contrast, \texttt{cdu} minimises raw driving time irrespective of
  waiting, which can increase idle periods at combined nodes.
  \texttt{vrd} is van-centric and leaves moped routes largely
  unconstrained beyond feasibility.
  The managerial choice among these four objectives depends on whether
  the operator prioritises shift completion time, fuel cost, or
  driver utilisation.
\end{remark}

\section{Graph-sparsification preprocessing}
\label{sec:preprocessing}

Before passing the model to the solver, the directed graph $G=(V,E)$ is
split into two vehicle-specific edge sets $E_{vn}$ and $E_m$, and
three preprocessing rules are applied.
All three rules are \emph{exact}: they remove only edges that cannot
appear in any feasible solution, so optimal solutions are preserved.
As a by-product, constraints C4 and C5 (vehicle-accessibility) become
redundant and can be dropped, reducing the constraint count
proportionally.

\subsection{Vehicle-accessibility edge removal}
\label{subsec:pre:access}

Every edge incident to at least one node that is inaccessible to the
corresponding vehicle type is removed.
Formally, for the van graph any edge $(u,v)\in E$ with
$u\notin V_{vn}$ or $v\notin V_{vn}$ is deleted from $E_{vn}$; the
symmetric rule applies to $E_m$ with respect to $V_m$.
This rule alone eliminates all edges to and from van-inaccessible
nodes (e.g.\ pedestrianised streets or narrow alleys reachable only
by mopeds), and vice versa.

\subsection{Time-window infeasibility edge removal}
\label{subsec:pre:tw}

A directed edge $(u,v)$ is \emph{time-window infeasible} for a given
vehicle type if the vehicle could never arrive at $v$ within its time
window after completing service at $u$:
\begin{equation}
  b_u + s_u + c_{(u,v)} > b_v,
  \label{tw:prune}
\end{equation}
where $c_{(u,v)}$ denotes the travel time for the vehicle in question.
Condition~\eqref{tw:prune} is sufficient: even if the vehicle arrived
at $u$ at the latest permissible time $b_u$ and departed immediately
after service, it would still miss the closing window at $v$.
The rule is applied independently to $E_{vn}$ and $E_m$ using their
respective travel-time matrices.

Note that combined nodes are exempt in the moped graph because
constraint~\eqref{c20} already forbids a moped from departing a
combined node before the van arrives; the synchronisation constraint
implicitly handles feasibility at these nodes.

The effectiveness of this rule depends heavily on the \emph{coverage
ratio}
\begin{equation}
  \varphi(\mathcal{W}) = \frac{s_\mathcal{W} \cdot b_\mathcal{W}}{t_\mathcal{W}},
  \label{phi}
\end{equation}
where $t_\mathcal{W}$, $s_\mathcal{W}$, and $b_\mathcal{W}$ are the
shift length, window width, and number of equally-spaced windows,
respectively.
Smaller $\varphi$ creates dead service zones between windows, making
more edges infeasible and therefore removable.

\subsection{Consecutive combined-node moped-edge removal}
\label{subsec:pre:consec}

If two nodes $u$ and $v$ both have demand exceeding the moped capacity
($d_u > L_m$ and $d_v > L_m$), then both must be served by the van and
are therefore combined nodes.
Any moped edge $(u,v)\in E_m$ between two combined nodes is removed,
because a moped arriving at $u$ would already have used its full
capacity on the goods loaded at $u$; it cannot carry anything further
to $v$ without first visiting at least one non-combined customer.
Formally:
\begin{equation}
  (u,v)\in E_m \text{ is removed if } d_u > L_m \text{ and } d_v > L_m.
  \label{consec:rule}
\end{equation}

\subsection{Reduction statistics}
\label{subsec:pre:stats}

The combined effect of the three rules was measured on instances
$\mathcal{I}_\mathrm{L}(30, \varphi, 4)$ drawn from the Lisbon network
(30 nodes, moped capacity $L_m=4$, coverage ratio $\varphi$), with
$|V_{vn}|/|V|$ and $|V_m|/|V|$ varied jointly.
Even when both ratios equal 1 (no delivery grey zones), graph
sparsification reduces the constraint count by between 12.5\%
and 25\% depending on $\varphi$.
The reduction grows as $\varphi$ decreases: at $\varphi=0.5$ the dense
dead zones between narrow time windows allow far more edge deletions.
Conversely, at $\varphi>1$ (overlapping windows) few edges are removed
by Rule~2, because visiting nodes in reverse temporal order is
feasible more often.
This confirms that graph sparsification is most effective precisely
when the model is hardest to solve (tight, non-overlapping time
windows), providing its greatest benefit where it is most needed.

\section{Cluster-based combined routing heuristic}
\label{sec:heuristic}

\subsection{Motivation and overview}
\label{subsec:heur:overview}

For instances with more than roughly 50 customers the exact MILP
becomes unreliable within practical time budgets (Section~\ref{subsec:exp:exact}).
We therefore develop a \emph{Cluster-Based Combined Routing} (CBCR)
heuristic inspired by the clustering approach of
\citet{ctbrp} for the CTBRP.
The key idea is to aggregate customers into $k$ macro-nodes (clusters),
solve the CVMRPTW exactly on the much smaller cluster graph, and then
expand the cluster-level routes back to individual customer routes
through a series of \emph{Shortest Hamiltonian Path} (SHP) subproblems.

CBCR proceeds in three phases:
\begin{enumerate}
  \item \textbf{Instantiation.} Customers are grouped into $k$ clusters
        and the instance parameters (time windows, service times, travel
        costs) are aggregated to cluster level.
        The MILP is solved on the cluster graph to obtain a combined
        cluster route $(\mathcal{V}', \mathcal{M}')$.
  \item \textbf{Van route reconstruction.} The actual customer-level
        van route is built by solving an SHP subproblem for each
        inter-combined-cluster fragment, enforcing that nodes inside
        combined clusters are visited last within each fragment.
  \item \textbf{Moped route reconstruction.} For each moped route
        segment in $\mathcal{M}'$, an SHP subproblem constructs the
        detailed customer-level moped path, linking its end point to the
        van arrival time at the next combined cluster to satisfy
        synchronisation.
\end{enumerate}

\subsection{Clustering and parameter aggregation}
\label{subsec:heur:cluster}

The number of clusters is set to
\begin{equation}
  k = \left\lceil \frac{\sum_{v\in V_c} d_v}{L_m} \right\rceil,
  \label{cbcr:k}
\end{equation}
so that each cluster has an expected demand close to $L_m$.
By default, $k$-means clustering (applied to customer coordinates)
is used.
Each cluster $\mathcal{C}_i$ is then mapped to a proxy node $v_i$
whose parameters are derived as follows:
\begin{itemize}
  \item \textbf{Travel costs.} $c^{vn}_{(v_i,v_j)} = \mathrm{avg}_{m\in\mathcal{C}_i,\,n\in\mathcal{C}_j} c^{vn}_{(m,n)}$; analogously for $c^m$.
  \item \textbf{Time window.} $a_i = \min_{v\in\mathcal{C}_i} a_v$,
        $b_i = \min_{v\in\mathcal{C}_i} b_v$ (the tightest window to
        avoid stranding any customer).
  \item \textbf{Service time and demand.} $s_i = \sum_{v\in\mathcal{C}_i} s_v$,
        $d_i = \sum_{v\in\mathcal{C}_i} d_v$.
  \item \textbf{Vehicle accessibility.} $v_i\in V_{vn}$ iff every node
        in $\mathcal{C}_i$ is van-accessible; similarly for $V_m$.
\end{itemize}

A cluster with $d_i > L_m$ is treated as a combined node in the
cluster-level solve.
Because all mopeds must pass through combined clusters, the fraction
of forced combined clusters grows as $k$ decreases (fewer, larger
clusters), making the cluster-level problem easier to solve but the
intra-cluster expansion more constrained.

To account for unknown intra-cluster travel, the time-propagation
constraints in the cluster-level model are augmented with additive
offsets $\texttt{min\_c}_c$ (van) and $\gamma\cdot\texttt{max\_c}_c$
(moped), where
\begin{equation}
  \texttt{min\_c}_c = \min_{s\in\mathcal{C}_c}
    \bigl\{ t_{\mathrm{end}} - a^s_v :
      \mathrm{SHP}(\mathcal{C}_c\setminus\{s\}, [(\_,\,s,\,a^s_v)]) \bigr\},
  \quad
  \texttt{max\_c}_c = \max_{s\in\mathcal{C}_c}(\cdots),
\end{equation}
represent the shortest and longest intra-cluster Hamiltonian paths
starting from each possible start node.
The parameter $\gamma\geq 1$ controls the conservatism of the moped
time estimate; setting $\gamma$ too large may make the cluster-level
problem infeasible, while $\gamma$ too small can lead to
synchronisation failures in the expansion phase.
Based on empirical testing, $\gamma=1.3$ provides a good trade-off.

\subsection{Route reconstruction via shortest Hamiltonian path}
\label{subsec:heur:reconstruct}

The SHP subroutine (Algorithm~\ref{alg:shortest-hamiltonian-path-paper})
solves a small TSP-like MILP to find the cheapest Hamiltonian path
through a given set of nodes \texttt{visit}, starting from a set of
candidate start points \texttt{starts} (each annotated with a cluster
index and an earliest-departure time), optionally forcing a specific
subset \texttt{last\_to\_visit} to be visited last, and optionally
ending at a set \texttt{ends}.
Its objective is to minimise the sum of arrival times:
$\min \sum_{v\in\texttt{visit}\cup\texttt{ends}} T_v$.
Time windows and van-synchronisation constraints are enforced locally.

\begin{algorithm}[!tbp]
\caption{ShortestHamiltonianPath (SHP) — sketch.}
\label{alg:shortest-hamiltonian-path-paper}
\small
\DontPrintSemicolon
\KwIn{\texttt{visit}, \texttt{starts}, \texttt{last\_to\_visit}, \texttt{ends}, $T^{vn}$, $c$, $a$, $b$, $s$}
\KwOut{\texttt{fragment}, $t_{\mathrm{start}}$, $t_{\mathrm{end}}$, \texttt{arrivals}}
Build a small MILP with: one selected start from \texttt{starts};
each node in \texttt{visit} visited exactly once; nodes in
\texttt{last\_to\_visit} visited after all others; if \texttt{ends}$\neq\emptyset$,
exactly one end node selected and synchronised to van arrival time.
  Solve with \texttt{grp0} (see Section~\ref{subsec:exp:setup}) and return optimal route and arrival times.
\end{algorithm}

For the van route, SHP is called once per inter-combined-cluster
fragment, with the nodes of the \emph{next} combined cluster as
\texttt{last\_to\_visit} so that the van arrives at the transfer point
after serving all intermediate customers.

For the moped routes, SHP is called once per moped segment
(the cluster sequence between two consecutive combined clusters in
$\mathcal{M}'$).
If a moped is \emph{linked} (i.e.\ it continues from a previous
segment), its last visited node is used as the sole start; otherwise,
any node in the starting combined cluster may be chosen
(\emph{fresh start}).
When $n$ mopeds start their shift from the same combined cluster, these $n$
mopeds always use a fresh start, while any remaining ones use linked starts.

\subsection{Complexity and parameter sensitivity}
\label{subsec:heur:complexity}

The dominant cost of CBCR is the collection of SHP calls, each of
which is NP-hard in general but involves only $O(|\mathcal{C}_i|)$
nodes per cluster.
In practice the SHP subproblems are solved in a few seconds with short
time budgets ($T_1=\SI{15}{\second}$, $T_2=\SI{3}{\second}$), making
the heuristic tractable for instances with up to 80--90 customers
(Section~\ref{subsec:exp:cbcr}).

The main sensitivity parameters are $k$ (cluster count), $\gamma$
(moped time buffer), and the clustering algorithm.
Smaller $k$ reduces the cluster-level solve time but amplifies
intra-cluster uncertainty; larger $k$ increases the cluster-level
problem size.
The choice $k=\lceil D/L_m\rceil$ (total demand divided by moped
capacity) keeps the expected cluster demand close to $L_m$, balancing
these two effects.

\subsection{Post-hoc validity}
\label{subsec:heur:a_posterior_validaty}

A solution is feasible only if it passes a post-hoc validity check given
in \eqref{def:validity}. Temporal and structural constraints are
verified. While the model is constructed to enforce all relevant constraints a
priori, this facilitates the development of further models and heuristics, such
as the CBCR heuristic. Let $\mathcal{P}(\mathcal{V}, \mathcal{M})$ denote the
set of conditions to be verified:
\begin{align}
\label{def:validity}
\begin{cases}
\mathcal{V}_1 = v_s,\;\; \mathcal{V}_{|\mathcal{V}|} = v_e, & \\[0.5ex]
\{\mathcal{V} \cup \bigcup_{i=1}^{|\mathcal{M}|} \mathcal{M}_i\} = \{V\}, \\[0.5ex]
|\{\mathcal{V}\}| = |\mathcal{V}|, & \\[0.5ex]
|\{f\}| = |f|, 
& \forall\, f \in \bigcup_{i=1}^{|\mathcal{M}|} \mathcal{Z}(\mathcal{M}_i), \\[0.5ex]
v \notin V \setminus V_{vn}, & \forall v \in \mathcal{V}, \\[0.5ex]
v \notin V \setminus V_m, & \forall v \in \bigcup_{i=1}^{|\mathcal{M}|} \mathcal{M}_i, \\[0.5ex]
|\mathcal{M}| = y^m, & \\[0.5ex]
\mathcal{M}_{i_1} \in K,\ \mathcal{M}_{i_{|\mathcal{M}_i|}} \notin K, & \forall \mathcal{M}_i \in \mathcal{M},\\[0.5ex]
\sum_{n \in f} d_n \le L_m, 
& \forall\, f \in \bigcup_{i=1}^{|\mathcal{M}|} \mathcal{Z}(\mathcal{M}_i), \\[0.5ex]
T^{vn}_u + s^{vn}_u + c^{vn}_{(u,v)} \;\le\; T^{vn}_v, & \forall (u,v) \in \mathcal{V},\ v \succ u,\\[0.5ex]
T^{vn}_v \;\in\; [\,a_v,\, b_v\,], & \forall v \in \mathcal{V},\\[0.5ex]
T^m_k \leq T^{vn}_k, & \forall k \in K, \\[0.5ex]
T^m_v \;\ge\; T^{vn}_u +
\begin{cases}
 s^{vn}_u + c^m_{(u,v)}, & u \in K,\\
 s^m_u + c^m_{(u,v)}, & u \notin K,
\end{cases} & \forall (u,v) \in \bigcup_{i=1}^{|\mathcal{M}|} \mathcal{M}_i,\ v \succ u, \\[1ex]
T^m_v \;\in\; [\,a_v,\, b_v\,], & \forall v \in \bigcup_{i=1}^{|\mathcal{M}|} \mathcal{M}_i \setminus K, \\[0.5ex]
\end{cases}
\end{align}
where $\mathcal{Z}$ fragments a moped route by removing its
combined nodes.

Interestingly, $\mathcal{P}$ also works for the common-depot variant of the
model. For the active-waiting-van variant we replace the 10th, 11th,
12th, and 13th conditions by
\begin{equation}
\label{def:validity2}
\left\{
\begin{array}{ll}
\max(T^{vn}_u, a_u) + s^{vn}_u + c^{vn}_{(u,v)} \le T^{vn}_v, 
& \forall (u,v) \in \mathcal{V},\ v \succ u,\\[0.3em]
T^{vn}_v \leq b_v, 
& \forall v \in \mathcal{V},\\[0.3em]
T^m_k \leq \max(T^{vn}_k, a_k), 
& \forall k \in K,\\[0.8em]
T^m_v \;\ge\; \max(T^{vn}_u, T^m_u) + 
\begin{cases}
 s^{vn}_u + c^m_{(u,v)}, & u \in K, \\[0.3em]
 s^m_u + c^m_{(u,v)}, & u \notin K,
\end{cases}
& \forall (u,v) \in \bigcup_{i=1}^{|\mathcal{M}|} \mathcal{M}_i,\ v \succ u.
\end{array}
\right.
\end{equation}

\section{Computational experiments}
\label{sec:experiments}

\subsection{Instance generation}
\label{subsec:exp:instances}

Instances are drawn from two real-world networks: \textbf{Lisbon}
(Portugal) and \textbf{Stuttgart} (Germany).
Both cities have comparable population sizes but markedly different
urban morphologies: Lisbon has a compact, irregular centre with narrow
streets, while Stuttgart has a more structured layout with wider
arterial roads.
Road networks and average speeds (differentiated by road type and
vehicle class) are obtained from OpenStreetMap \citep{openstreetmap}.

An instance $\mathcal{I}_\mathrm{L}(n, \varphi, L_m)$ is generated as
follows.
A city centre is fixed and a radius limiting the maximum one-way travel
time from the centre is set to \SI{30}{\minute}.
$n-2$ customer locations are sampled uniformly within this area.
Time windows are drawn from
$\mathcal{W}(180, \varphi\cdot 180/3, 3)$ (shift of 180~min,
three equally-spaced windows of width $\varphi\cdot 60$~min),
assigned with uniform probability.
Customer demand is drawn uniformly from $\{1,\dots,L_m+1\}$, so that
some customers require van service.
The subscripts $\mathrm{L}$ and $\mathrm{S}$ distinguish Lisbon and
Stuttgart instances.
To guard against infeasibility caused by tight early windows at distant
nodes, the generator iteratively relaxes demand constraints on the
hardest-to-reach forced van-stops until the minimal van route is
feasible.

Standard solver settings use $f_1=\min y^m$,
$f_2=\texttt{tct}$, first-stage budget $T_1=\SI{600}{\second}$,
second-stage budget $T_2=\SI{30}{\second}$, MIP gap $0.02$ for
$n\leq30$, $0.10$ for $31\leq n\leq50$, and $0.20$ for $n>50$;
the MIP gap determines the minimal quality of the solution, with a higher value yielding a worse solution faster.
Each configuration is run on 10 independent instances.
Solutions were validated according to Section~\ref{subsec:heur:a_posterior_validaty}.

\subsection{Solver configuration and hardware}
\label{subsec:exp:setup}

All experiments were conducted on a virtual machine (KVM,
Debian GNU/Linux 12) equipped with an Intel\textsuperscript{\textregistered}
Xeon\textsuperscript{\textregistered} Silver 4214 at \SI{2.20}{\giga\hertz}
(10 logical cores), \SI{26}{\gibi\byte} RAM, and \SI{80}{\mebi\byte} L3
cache.
Gurobi 12.0.2~\citep{gurobi} was used with Python~3.11.2.
To allow triple-test parallelisation, the thread count was set to 3 and
\texttt{NodefileStart} to 6~GB.

Four Gurobi configurations are compared (Table~\ref{tab:gurobi}).
\texttt{grp0} is conservative (limited heuristics, no cuts);
\texttt{grp1} is balanced;
\texttt{grp2} is aggressive (high heuristics, maximum cuts);
\texttt{grp3} prioritises primal feasibility early (\texttt{MIPFocus}~2).

\begin{table}[!tbp]
  \centering
  \caption{Gurobi search configurations and optimality settings.}
  \label{tab:gurobi}
    \centering
    \begin{tabular}{ccccc}
      \toprule
      Config. & MIPFocus & Heuristics & Presolve & Cuts \\
      \midrule
      \texttt{grp0} & 1 & 0.2 & 1 & 0 \\
      \texttt{grp1} & 1 & 0.3 & 2 & 1 \\
      \texttt{grp2} & 1 & 0.9 & 2 & 3 \\
      \texttt{grp3} & 2 & 0.4 & 1 & 1 \\
      \bottomrule
    \end{tabular}
    \vskip1ex
    \begin{tabular}{ccc}
      \toprule
      $n$ & MIP gap & NodefileStart \\
      \midrule
      $\leq 30$ & 0.02 & 6 GB \\
      $31$--$50$ & 0.10 & 6 GB \\
      $> 50$ & 0.20 & 6 GB \\
      \bottomrule
    \end{tabular}
\end{table}

\subsection{Performance metrics}
\label{subsec:exp:metrics}

Beyond $y^m$ and $f_2$, three solution-quality metrics are reported.

\medskip
\noindent\textbf{Average moped load (AML).}
The mean fraction of moped capacity used across all inter-combined-node
route segments:
\begin{equation}
  \mathrm{AML} = \mathrm{avg}_{f\in\mathcal{Z}(\bigcup_i \mathcal{M}_i)}
    \!\left(\frac{\sum_{v\in f} d_v}{L_m}\right),
\end{equation}
where $\mathcal{Z}(r)$ fragments a moped route $r$ by removing its
combined nodes.

\medskip
\noindent\textbf{Average synchronisation ratio (ASR).}
The fraction of total active time spent waiting or synchronising rather
than travelling:
\begin{equation}
  \mathrm{ASR} =
  \frac{(T_{v_e}-T_{v_s} - \sum_{(i,j)\in\mathcal{V}}c^{vn}_{(i,j)})
        + \sum_i(T^b_{\mathcal{M}_{i,|\mathcal{M}_i|}} - T^b_{\mathcal{M}_{i,1}}
          - \sum_{(i,j)\in\mathcal{M}_i}c^m_{(i,j)})}
  {(T_{v_e}-T_{v_s})
   + \sum_i(T^b_{\mathcal{M}_{i,|\mathcal{M}_i|}} - T^b_{\mathcal{M}_{i,1}})}.
\end{equation}

\medskip
\noindent\textbf{Active ratios.}
The van active ratio $A_\mathcal{V}$ and moped active ratio
$A_\mathcal{M}$ decompose ASR into per-fleet driving fractions.

\subsection{Exact model: scalability and solution quality}
\label{subsec:exp:exact}

Table~\ref{tab:gurobi_results_quality} summarises solution quality for
increasing $n$ on Lisbon instances ($\varphi=1$, $L_m=4$) under
\texttt{grp0}, which outperforms the other three configurations.
The aggressive \texttt{grp2} and moderate \texttt{grp1} performed
worst; \texttt{grp3} matched \texttt{grp0} on solve rate but produced
slightly higher $y^m$.
This outcome is consistent with the observation that focusing on
the dual bound (lower bound) via cuts tends to be less effective than
the conservative primal-first strategy for this class of problem.

The moped count grows approximately linearly with $n$ for $n>18$,
which is operationally desirable: splitting a shift asymmetrically
across two sub-shifts costs no more than combining them.
Optimal proofs are lost beyond $n=13$; the solve rate drops below
80\% at $n=48$ and becomes unreliable beyond that.
AML saturates at 0.9 and ASR at 0.36, indicating that larger instances
make fuller use of each moped's capacity and achieve more consistent
synchronisation.

\begin{table}[!tbp]
\centering
\caption{Solution quality for $\mathcal{I}_\mathrm{L}(n,1,4)$ under
         \texttt{grp0}: mean mopeds ($y^m$), secondary objective
         ($f_2 = T(\mathcal{V})+T(\mathcal{M})$), AML, ASR, and
         percentage of optimally and feasibly solved instances.}
\label{tab:gurobi_results_quality}
\resizebox{0.62\linewidth}{!}{%
\begin{tabular}{ccccccc}
\toprule
$n$ & $y^m$ & $f_2$ & AML & ASR & Optimal & Solved \\
\midrule
 8 & 0.1 & 136.0 & n.s. & 0.65 & 100 & 100 \\
13 & 0.1 & 160.4 & n.s. & 0.49 & 100 & 100 \\
18 & 1.1 & 192.0 & 0.77 & 0.42 &   0 & 100 \\
23 & 1.2 & 244.5 & 0.82 & 0.38 &   0 & 100 \\
28 & 1.9 & 318.8 & 0.79 & 0.38 &   0 & 100 \\
33 & 2.3 & 414.4 & 0.87 & 0.37 &   0 &  90 \\
38 & 3.2 & 505.8 & 0.87 & 0.38 &   0 &  90 \\
43 & 4.2 & 564.4 & 0.90 & 0.36 &   0 &  80 \\
48 & 5.0 & 676.8 & 0.89 & 0.36 &   0 &  80 \\
\bottomrule
\end{tabular}}
\end{table}

Model variables and non-zeros grow as $O(n^2)$ (due to edge variables),
while constraints grow as $O(n)$.
This quadratic growth drives the observed memory-use and runtime
escalation: the first-stage budget of 600~s is exhausted for all
$n\geq23$, and RAM consumption reaches over 500~MB for $n=48$.

\subsection{Effect of instance parameters}
\label{subsec:exp:params}

\begin{figure}[!tbp]
  \centering
  \begin{subfigure}[t]{0.5\textwidth}
    \centering
    \includegraphics[width=\linewidth]{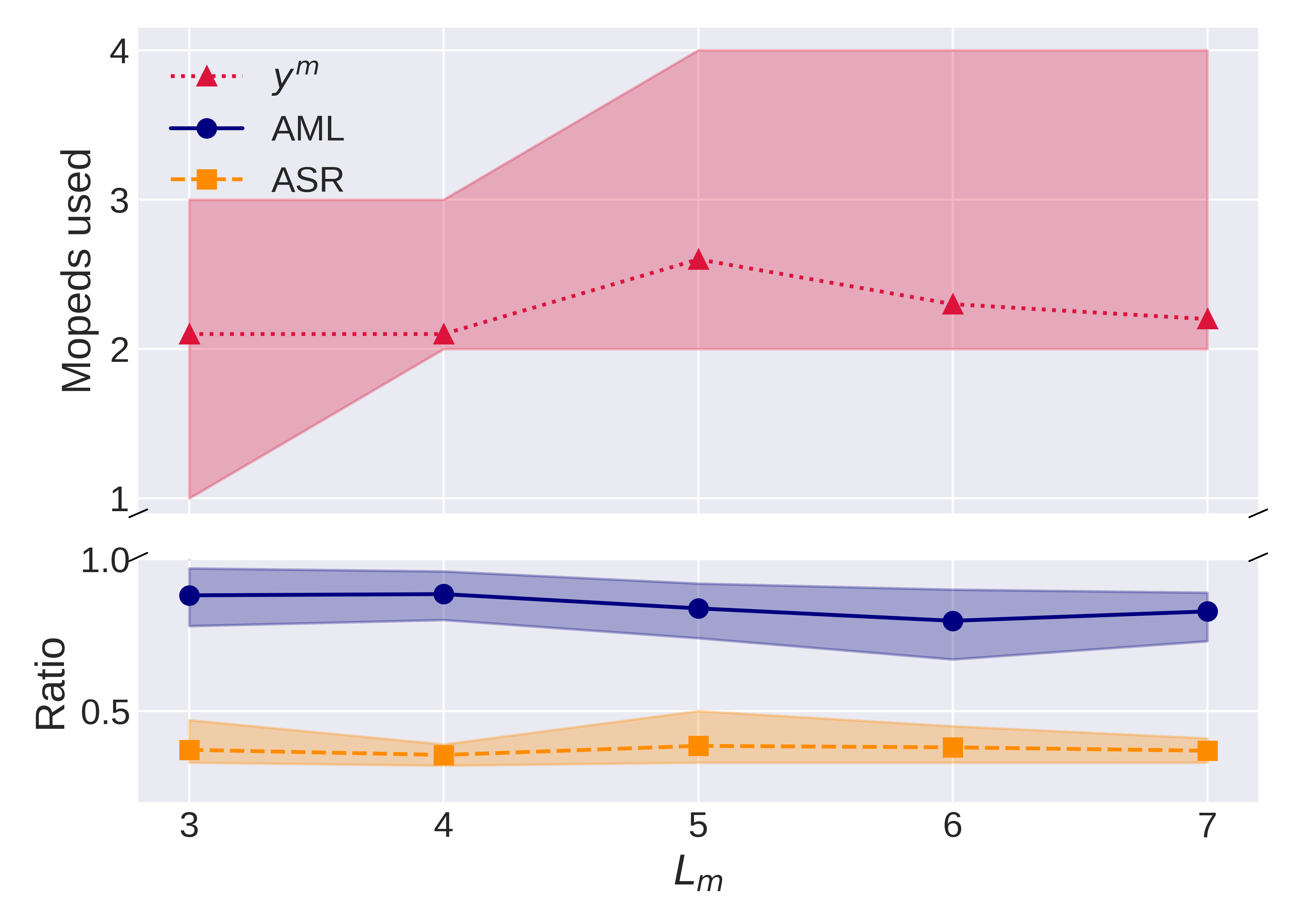}
    \caption{$\mathcal{I_L}(30, 1, L_m)$}
    \label{fig:moped_capacity}
  \end{subfigure}\hfill
  \begin{subfigure}[t]{0.5\textwidth}
    \centering
    \includegraphics[width=\linewidth]{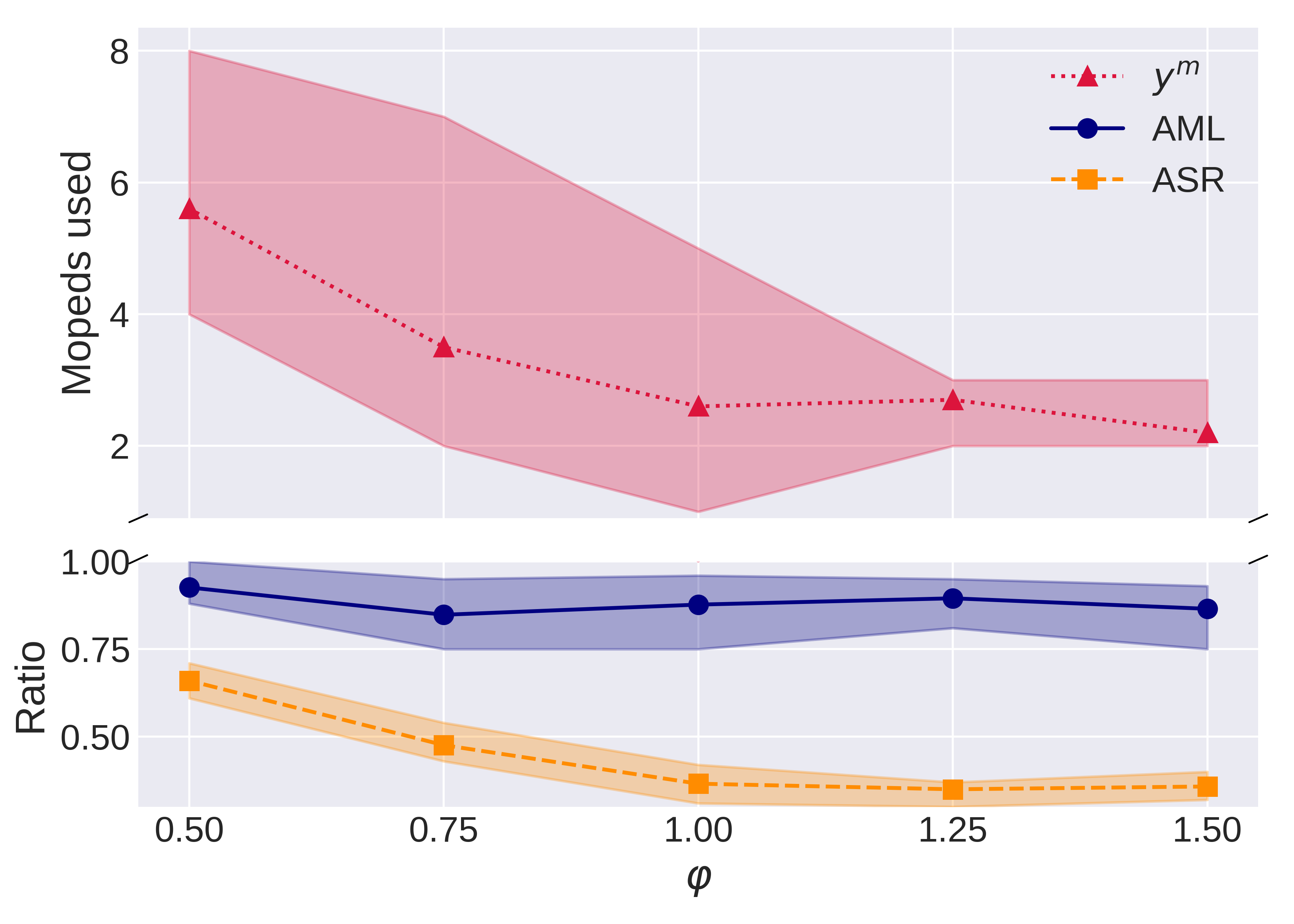}
    \caption{$\mathcal{I_L}(30, \varphi, 4)$}
    \label{fig:moped_time_windows}
  \end{subfigure}

  \caption{Mean and range of the number of mopeds used (\(y^m = f_1\)), average moped load (AML), and average synchronization ratio (ASR) for varying moped capacities and coverage ratios.}
  \label{fig:moped_capacity_moped_time_windows}
\end{figure}

\paragraph{Moped capacity $L_m$}
Increasing $L_m$ from 3 to 7 does not reduce $y^m$ substantially,
because demand is drawn from $\{1,\dots,L_m+1\}$, keeping the ratio
$d/L_m$ roughly constant (see Figure~\ref{fig:moped_capacity}).
AML decreases with $L_m$ (larger capacity is harder to fill
completely), while ASR remains unaffected.

\paragraph{Coverage ratio $\varphi$}
As $\varphi$ falls below 1.0, dead service zones between non-overlapping
windows force more conflicting customer pairs to use separate mopeds,
causing $y^m$ to increase sharply and ASR to nearly double (see Figure~\ref{fig:moped_time_windows}).
AML does not depend significantly on $\varphi$.

\subsection{Problem variant comparison}
\label{subsec:exp:variants}

\begin{figure}[!hbp]
  \centering

  \begin{subfigure}[t]{0.48\textwidth}
    \centering
    \includegraphics[width=\linewidth]{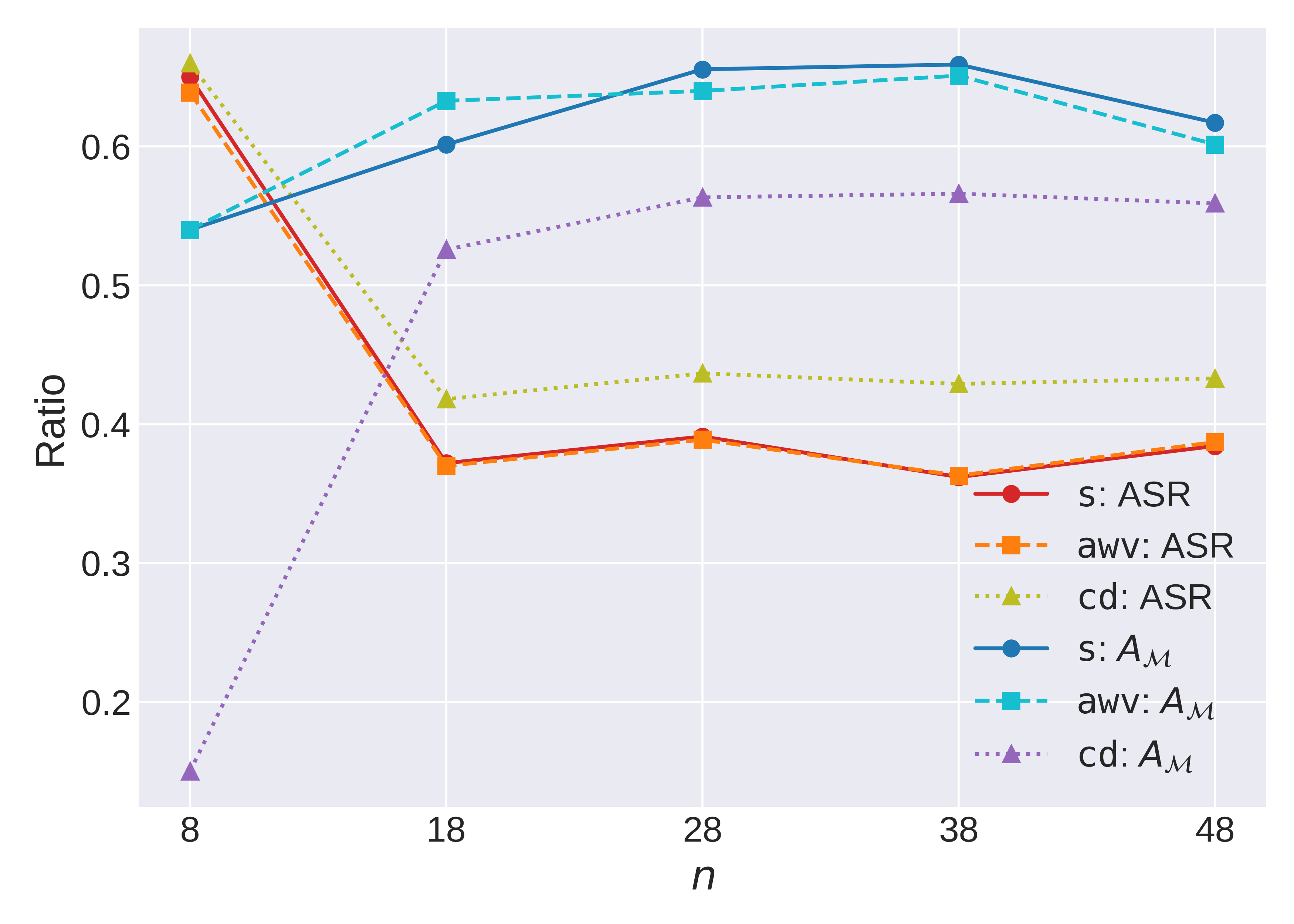}
    \caption{$\mathcal{I_L}(n, 1, 4)$}
    \label{fig:var11}
  \end{subfigure}\hfill
  \begin{subfigure}[t]{0.48\textwidth}
    \centering
    \includegraphics[width=\linewidth]{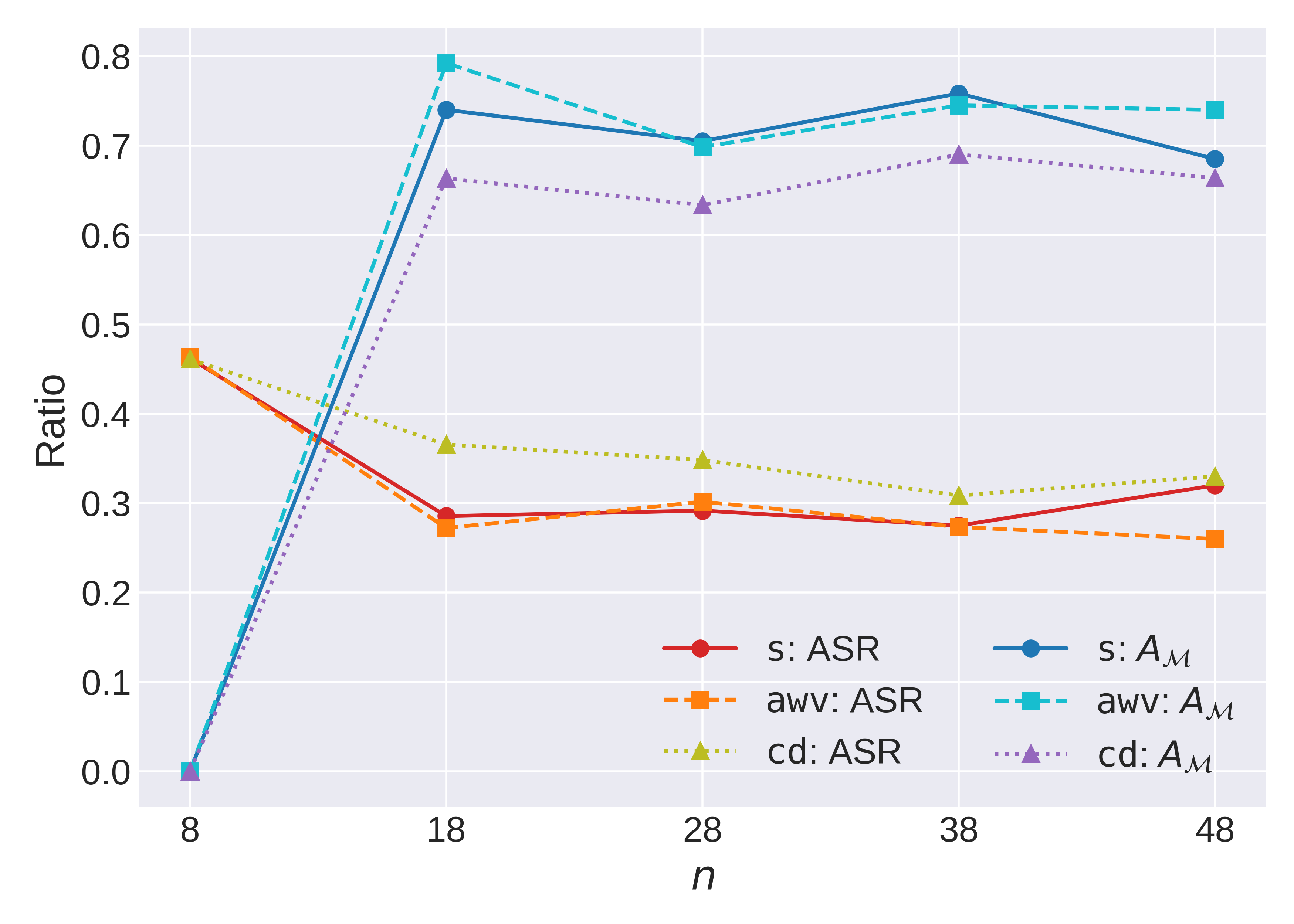}
    \caption{$\mathcal{I_S}(n, 1, 4)$}
    \label{fig:var12}
  \end{subfigure}

  \begin{subfigure}[t]{0.48\textwidth}
    \centering
    \includegraphics[width=\linewidth]{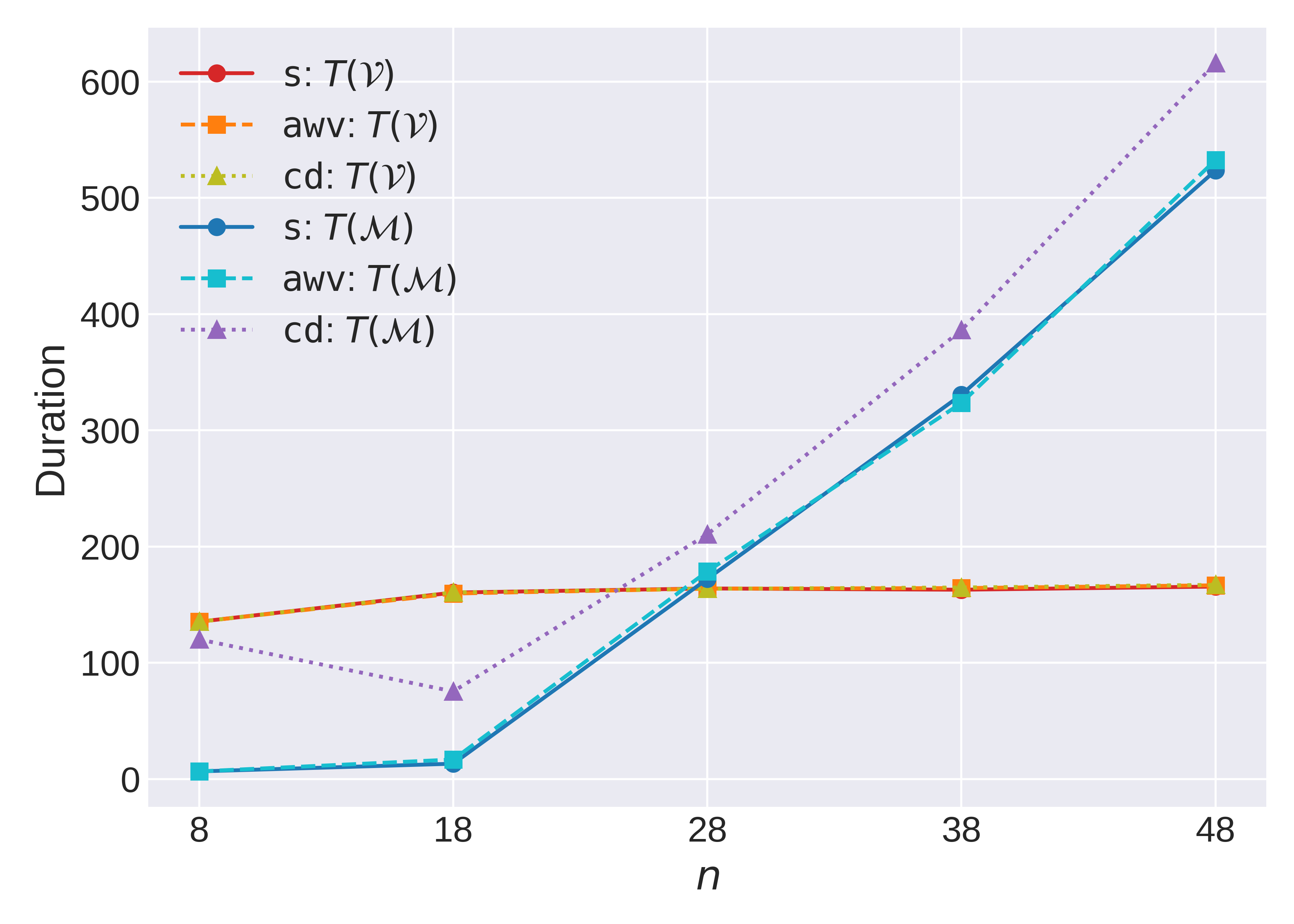}
    \caption{$\mathcal{I_L}(n, 1, 4)$}
    \label{fig:var21}
  \end{subfigure}\hfill
  \begin{subfigure}[t]{0.48\textwidth}
    \centering
    \includegraphics[width=\linewidth]{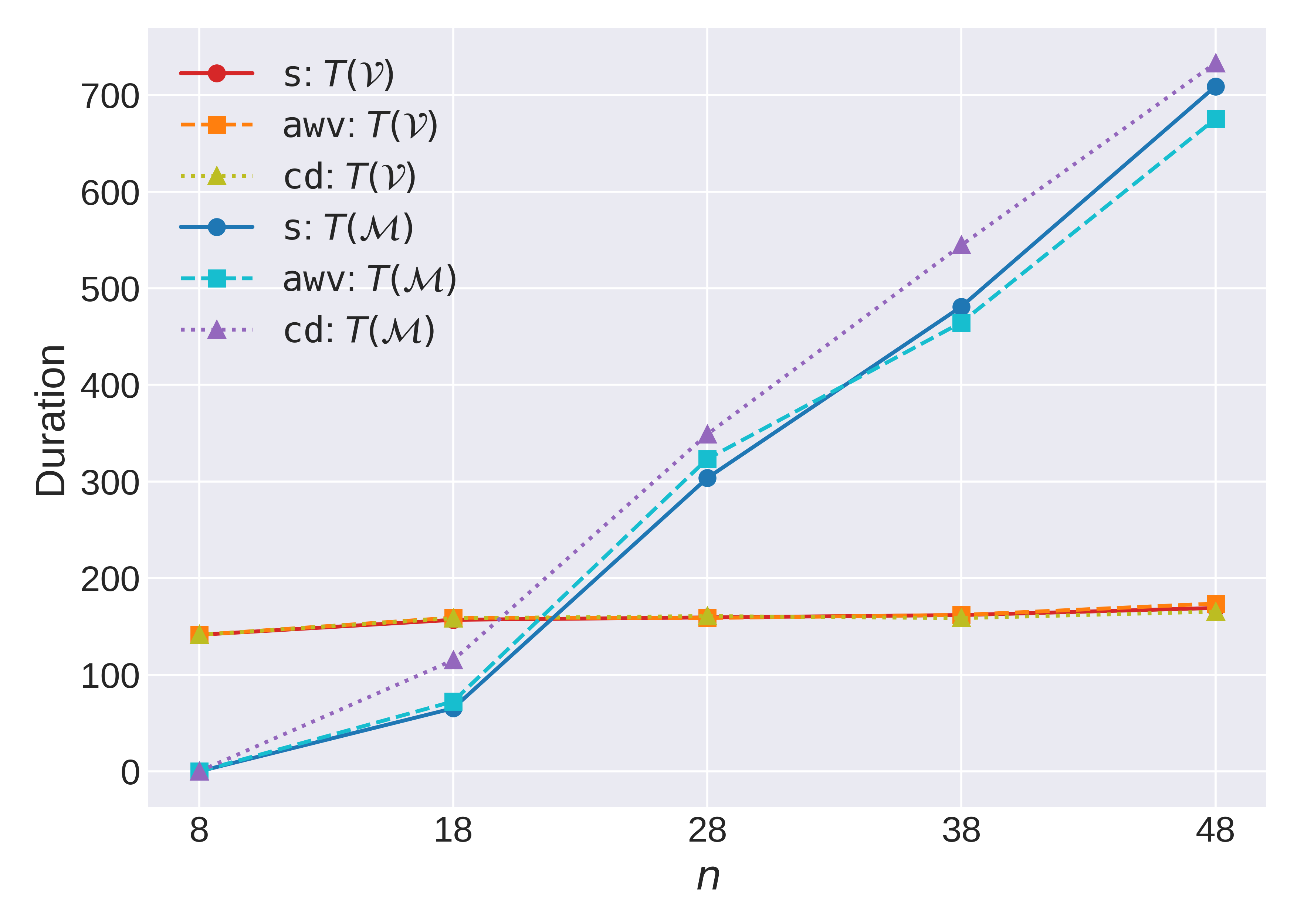}
    \caption{$\mathcal{I_S}(n, 1, 4)$}
    \label{fig:var22}
  \end{subfigure}

  \begin{subfigure}[t]{0.48\textwidth}
    \centering
    \includegraphics[width=\linewidth]{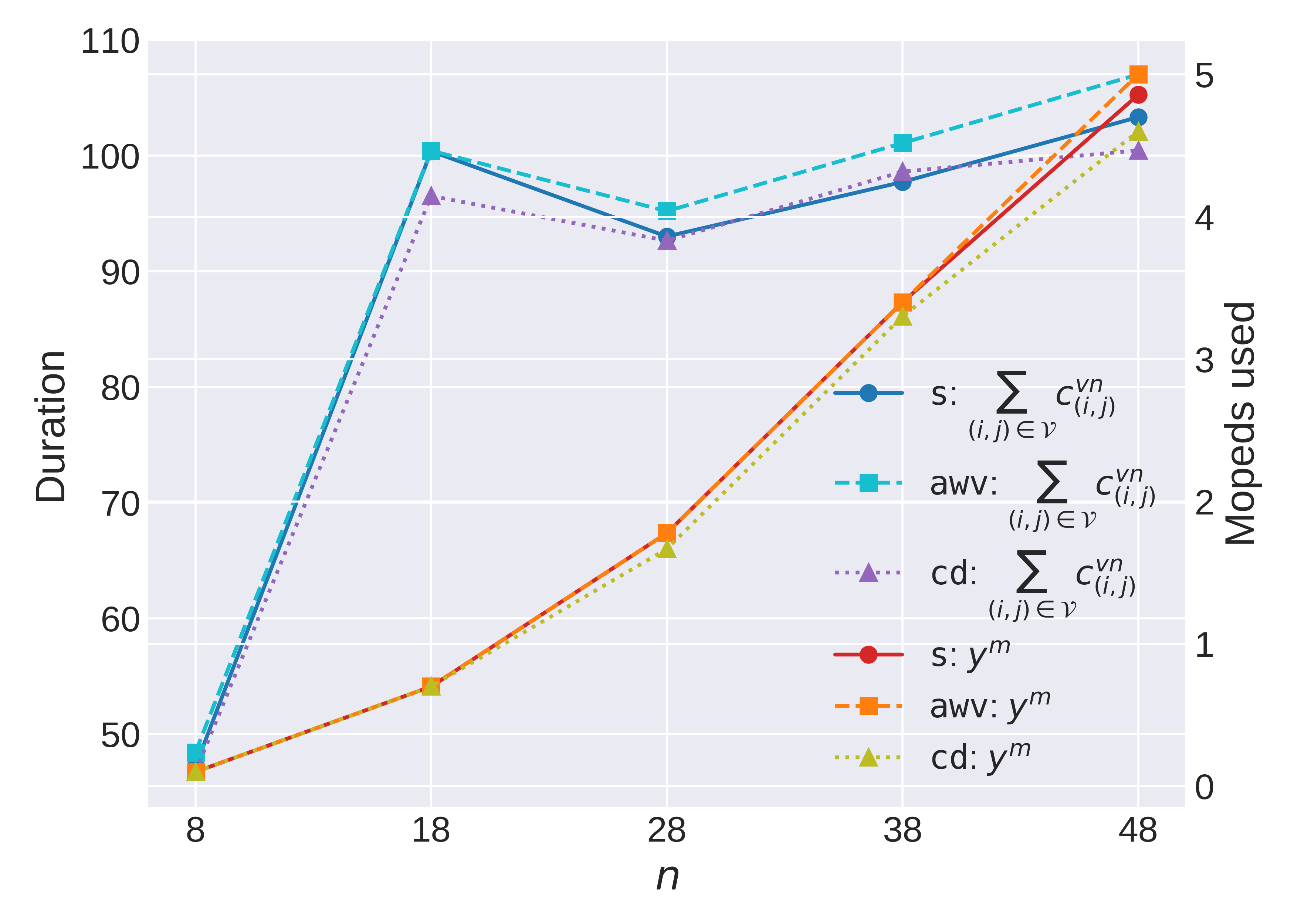}
    \caption{$\mathcal{I_L}(n, 1, 4)$}
    \label{fig:var31}
  \end{subfigure}\hfill
  \begin{subfigure}[t]{0.48\textwidth}
    \centering
    \includegraphics[width=\linewidth]{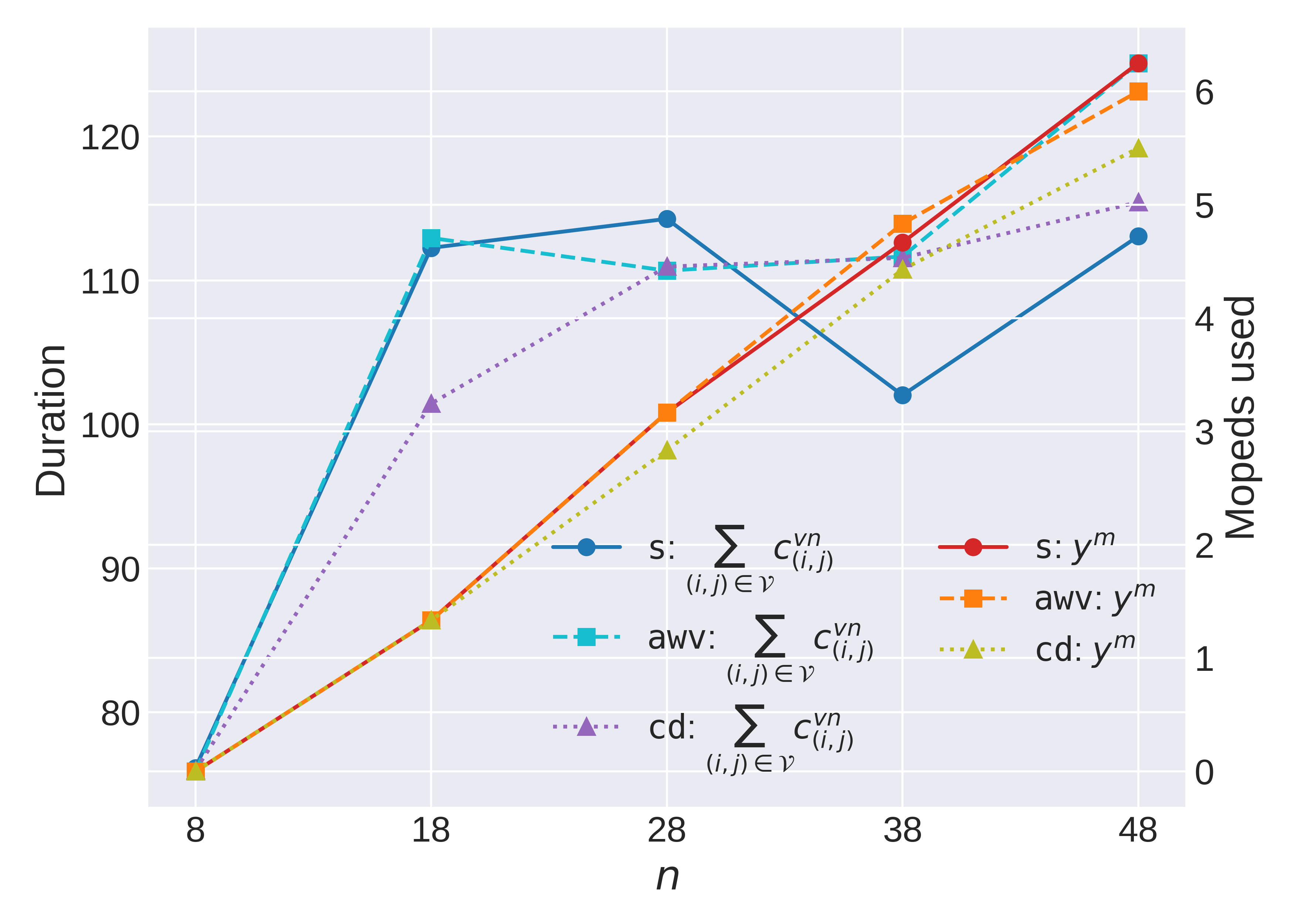}
    \caption{$\mathcal{I_S}(n, 1, 4)$}
    \label{fig:var32}
  \end{subfigure}

  \caption{Comparison of the standard, active-waiting-van, and common-depot variants:  
  (a,b) Mean average synchronization ratio (ASR) and moped active ratio ($A_\mathcal{M}$),  
  (c,d) Mean van and moped route duration (in~\si{\minute}),  
  (e,f) Mean van driving time (in~\si{\minute}) and number of mopeds used (\(y^m = f_1\)).}
  \label{fig:var_all}
\end{figure}

The three variants (standard \textsc{s}, active-waiting-van
\textsc{awv}, and common-depot \textsc{cd}) were compared on
$\mathcal{I}_\mathrm{L}(n,1,4)$ and $\mathcal{I}_\mathrm{S}(n,1,4)$ (see Figure~\ref{fig:var_all}).

\textsc{cd} achieves a slightly lower $y^m$ in some instances (all
mopeds can start earlier from the depot), but consistently shows
higher ASR and worse moped active ratio than \textsc{s} and
\textsc{awv}, because forcing all mopeds to depart from a single depot
increases their total idle time.
\textsc{awv} provides a modest relaxation over \textsc{s} by
decoupling the reload event from the van's delivery, allowing mopeds a
head-start.
However, the improvement is inconsistent: for the largest Stuttgart
instance only a slight gain is observed, and across most configurations
\textsc{s} and \textsc{awv} are largely equivalent.

Activity ratios are consistently higher in Stuttgart (more structured
road network), while $y^m$ and moped route durations are slightly
lower in Lisbon (more compact geography).

\subsection{Secondary objective comparison}
\label{subsec:exp:secondary}

Table~\ref{tab:lexicographic_improvement} reports the improvement
$\Delta f_2$ achieved in the second lexicographic stage relative to
the value at the end of the first stage, alongside the resulting van
route duration $T(\mathcal{V})$ and moped route duration
$T(\mathcal{M})$.

\texttt{tct} consistently delivers the largest improvement for the
mopeds ($\Delta T(\mathcal{M})$ up to $-75\%$) at the cost of a
slight increase in $T(\mathcal{V})$, because the solver marginally
lengthens the van route to create better moped synchronisation
windows.
Importantly, improving moped routes does not degrade the van route
significantly.
In Stuttgart, the \texttt{tct} objective achieves a moped active ratio
of almost 80\%.

\texttt{cdu} reduces total driving time by up to 37.4\% for $n=48$
in Lisbon, but leaves more waiting time in the schedule, resulting in
worse ASR and lower active ratios than \texttt{tct}.
\texttt{vrd} is van-centric and leaves moped routes largely
unconstrained; $A_\mathcal{M}$ under \texttt{vrd} is markedly lower
than under \texttt{tct}.

These results suggest that \texttt{tct} is the recommended secondary
objective for operators whose priority is fleet utilisation and shift
completion time.

\begin{table}[!tb]
\centering
\caption{Mean $\Delta f_2$ (percent improvement), van route duration
         $T(\mathcal{V})$, and moped route duration $T(\mathcal{M})$
         (both in minutes) for \texttt{vrd}, \texttt{tct}, and \texttt{cdu}
         on Lisbon and Stuttgart instances.}
\label{tab:lexicographic_improvement}
\small
\begin{tabular}{c ccc|ccc|ccc}
\toprule
\multicolumn{10}{c}{$\mathcal{I}_\mathrm{L}(n,1,4)$} \\[0.2em]
$n$
& \multicolumn{3}{c|}{\texttt{vrd}}
& \multicolumn{3}{c|}{\texttt{tct}}
& \multicolumn{3}{c}{\texttt{cdu}} \\
& $\Delta f_2$ & $T(\mathcal{V})$ & $T(\mathcal{M})$
& $\Delta f_2$ & $T(\mathcal{V})$ & $T(\mathcal{M})$
& $\Delta f_2$ & $T(\mathcal{V})$ & $T(\mathcal{M})$ \\
\midrule
 8 & 18.4 & 135.3 & 120.0 & 22.3 & 135.3 &   6.5 & 12.9 & 166.5 & 180.0 \\
18 &  8.7 & 150.9 & 162.8 & 33.1 & 160.3 &  13.1 & 29.4 & 166.8 & 105.2 \\
28 &  7.8 & 156.0 & 273.2 & 20.7 & 163.8 & 172.4 & 41.5 & 168.0 & 223.1 \\
38 &  8.4 & 154.5 & 514.1 & 25.7 & 162.7 & 332.2 & 51.1 & 168.4 & 471.6 \\
48 &  6.1 & 158.3 & 790.5 & 23.4 & 165.6 & 519.1 & 37.4 & 168.2 & 690.5 \\
\bottomrule
\end{tabular}
\vskip1ex
\begin{tabular}{c ccc|ccc|ccc}
\toprule
\multicolumn{10}{c}{$\mathcal{I}_\mathrm{S}(n,1,4)$} \\[0.2em]
$n$
& \multicolumn{3}{c|}{\texttt{vrd}}
& \multicolumn{3}{c|}{\texttt{tct}}
& \multicolumn{3}{c}{\texttt{cdu}} \\
& $\Delta f_2$ & $T(\mathcal{V})$ & $T(\mathcal{M})$
& $\Delta f_2$ & $T(\mathcal{V})$ & $T(\mathcal{M})$
& $\Delta f_2$ & $T(\mathcal{V})$ & $T(\mathcal{M})$ \\
\midrule
 8 & 12.0 & 141.4 &   0.0 & 12.0 & 141.4 &   0.0 &  9.5 & 162.8 &   0.0 \\
18 & 11.3 & 145.6 & 191.8 & 36.2 & 156.5 &  65.2 & 33.1 & 161.3 & 145.5 \\
28 &  7.7 & 149.6 & 485.5 & 28.6 & 159.4 & 303.6 & 40.4 & 161.9 & 407.9 \\
38 &  7.2 & 153.1 & 681.5 & 23.3 & 162.0 & 480.1 & 38.7 & 164.9 & 634.6 \\
48 &  6.8 & 158.6 & 958.9 & 23.8 & 169.2 & 703.6 & 30.1 & 169.8 & 833.2 \\
\bottomrule
\end{tabular}
\end{table}

\subsection{City comparison: Lisbon vs.\ Stuttgart}
\label{subsec:exp:cities}

Activity ratios are consistently higher in Stuttgart across all
variants and secondary objectives.
This reflects Stuttgart's wider, more regular road network, which
allows vehicles to travel more efficiently with less waiting at
combined nodes.
Total costs (moped count and route duration) are, however, marginally
lower in Lisbon, where the compact geography shortens inter-node
distances.
The gap in ASR between the two cities narrows as $n$ increases,
suggesting that the synchronisation overhead is primarily driven by
the number of combined interactions rather than the road structure.

\subsection{CBCR heuristic evaluation}
\label{subsec:exp:cbcr}

The CBCR heuristic was evaluated on a synthetic grid-map generator
(diagonal travel time $D_t=\SI{30}{\minute}$, a single 180-min window,
van/moped speed ratio 1.25, $\gamma=1.3$, $T_1=\SI{15}{\second}$,
$T_2=\SI{3}{\second}$, \textsc{awv} variant).
Table~\ref{tab:cbcr_tables} reports results for two demand
configurations.

With $L_m=6$ and demand $\sim\mathrm{Uniform}\{1,2,3\}$ (generous
capacity margin), CBCR solves 80\% of instances at $n=70$,
and 40\% at $n=80$ and 90, with AML around 0.85 and moderate ASR.
With $L_m=4$ and demand $\sim\mathrm{Uniform}\{1,2,3,4\}$ (tight
capacity), the solve rate drops to 20\% at $n=70$ and 80, and to 0\% at
$n=90$, and $y^m$ inflates dramatically (up to 13 mopeds at
$n=80$).

The primary failure mode is synchronisation infeasibility in the
expansion phase: as more clusters become forced combined nodes, the
intra-cluster path is over-constrained and no feasible SHP solution
exists.
Increasing $T_1$ and $T_2$ does not resolve this, as the bottleneck
is structural rather than computational.
Since many of the generated solutions were almost feasible and only slightly
violated timing constraints on a few edges, an improvement should be possible:
a per-cluster adaptive $\gamma$ strategy based on prior analysis of the graph
is identified as one possible direction for improvement.

\begin{table}[!tbp]
\centering
\caption{CBCR evaluation: mean $y^m$, AML, ASR, and solve rate (\%).}
\label{tab:cbcr_tables}
\begin{minipage}{0.47\linewidth}
  \centering
  \begin{tabular}{ccccc}
    \toprule
    $n$ & $y^m$ & AML & ASR & Solved \\
    \midrule
    $70$ & $2.2$ & $0.86$ & $0.49$ & $80.0$ \\
    $80$ & $3.0$ & $0.84$ & $0.54$ & $40.0$ \\
    $90$ & $5.0$ & $0.86$ & $0.57$ & $40.0$ \\
    \bottomrule
  \end{tabular}\\[1ex]
  $L_m=6$, $d\sim\mathrm{U}\{1,2,3\}$.
\end{minipage}
\begin{minipage}{0.47\linewidth}
  \centering
  \begin{tabular}{ccccc}
    \toprule
    $n$ & $y^m$ & AML & ASR & Solved \\
    \midrule
    $70$ & $7.0$ & $0.89$ & $0.75$ & $20.0$ \\
    $80$ & $13.0$ & $0.77$ & $0.61$ & $20.0$ \\
    $90$ & $0.0$ & $0.00$ & $0.00$ & $0.0$ \\
    \bottomrule
  \end{tabular}\\[1ex]
  $L_m=4$, $d\sim\mathrm{U}\{1,2,3,4\}$.
\end{minipage}
\end{table}

\section{Conclusions}
\label{sec:conclusions}

This paper introduced the Combined Van-and-Mopeds Routing Problem with
Time Windows (CVMRPTW), a new variant of the heterogeneous-fleet VRP
motivated by last-mile parcel delivery with cargo-moped assistance.
We developed a compact MILP formulation, discussed its NP-hardness, and
proposed three exact preprocessing rules that reduce the constraint
count by 12.5--25\% at no loss of optimality.
Two operational variants (\textsc{awv} and \textsc{cd}) extend the
base model, and a lexicographic two-stage framework supports four
secondary objectives.
A Cluster-Based Combined Routing heuristic extends tractability to
instances with up to 80 customers under favourable conditions.

The main empirical findings are:
\begin{enumerate}
  \item The moped count grows linearly with the number of customers,
        confirming that splitting a delivery shift across two sub-shifts
        is not penalised.
  \item The conservative Gurobi configuration (\texttt{grp0}) outperforms
        more aggressive settings, consistent with the observation that
        primal-first strategies are more effective than dual-bound-focused
        ones for this problem class.
  \item The \texttt{tct} secondary objective (minimise total completion
        time) is superior to \texttt{vrd} and \texttt{cdu} for fleet
        utilisation: it substantially reduces moped idle time without
        significantly lengthening the van route.
  \item The standard and active-waiting-van variants yield similar
        performance; the common-depot variant reduces the moped count
        slightly but incurs higher synchronisation overhead.
  \item City structure matters: Stuttgart's wider road network yields
        consistently higher active ratios, while Lisbon's compact
        geography produces lower moped counts and shorter routes.
  \item The CBCR heuristic is effective when the demand-to-capacity
        ratio is relaxed, but degrades when capacity is tight, pointing
        to the need for per-cluster adaptive time buffers.
\end{enumerate}

\paragraph{Managerial implications}
Urban logistics operators using van-and-moped fleets should adopt the
\texttt{tct} lexicographic objective to maximise vehicle utilisation
during a shift.
The standard synchronisation model is sufficient in most operational
settings; the \textsc{awv} relaxation provides only marginal gains
unless time windows are especially tight.
Graph sparsification should always be applied as a lightweight, exact
preprocessing step.
The CBCR heuristic is recommended for capacity-generous instances with
more than 50 customers; for tighter scenarios, more sophisticated
metaheuristics remain an open research direction.

\paragraph{Limitations and future work}
The model currently handles a single van; extending it to multiple vans
would require modifying the depot-flow constraints and re-assessing
the synchronisation logic.
The CBCR heuristic's performance is sensitive to $\gamma$ and the
cluster count $k$; an adaptive selection mechanism based on instance
structure is a natural next step.
Developing metaheuristics (e.g.\ large neighbourhood search or
adaptive large neighbourhood search) tailored to the CVMRPTW structure
would further extend scalability.
Finally, benchmarking on standardised last-mile delivery datasets would
allow direct comparison with existing heterogeneous-fleet VRP methods.

\section*{Acknowledgements}

The authors thank Dingoo for providing access to operational data and
for framing the real-world delivery problem that motivated this
research.
This work was supported in part by national funds through FCT ---
Funda\c{c}\~ao para a Ci\^encia e a Tecnologia --- under
projects UID/50021/2025 and UID/PRR/50021/2025 (INESC-ID).

\section*{CRediT author statement}

\textbf{Pedro Lameiras:} Conceptualization, Methodology, Software,
Formal analysis, Investigation, Writing -- original draft.
\textbf{Alexandre Francisco:} Conceptualization, Computational resources, Supervision, Writing -- review and editing.
\textbf{Adriano Serrano:} Conceptualization, Resources, Writing -- review and editing.
\textbf{C\'atia Vaz:} Formal analysis, Writing -- original draft, review and editing.

\section*{Declaration of competing interests}

The authors declare that they have no known competing financial interests
or personal relationships that could have appeared to influence the work
reported in this paper.

\bibliographystyle{elsarticle-harv}
\bibliography{./Bibliography}

@article{Miller1960,
  author    = {C. E. Miller and A. W. Tucker and R. A. Zemlin},
  title     = {Integer Programming Formulation of Traveling Salesman Problems},
  journal   = {Journal of the ACM},
  volume    = {7},
  number    = {4},
  year      = {1960},
  pages     = {326 -- 329},
  doi       = {10.1145/321043.321046},
  publisher = {ACM},
  address   = {New York, NY, USA},
}

@inbook{Karp1972,
  author    = {Richard M. Karp},
  title     = {Reducibility among Combinatorial Problems},
  journal   = {Complexity of Computer Computations},
  bookTitle = {Complexity of Computer Computations: Proceedings of a symposium on the Complexity of Computer Computations, held March 20--22, 1972, at the IBM Thomas J. Watson Research Center, Yorktown Heights, New York, and sponsored by the Office of Naval Research, Mathematics Program, IBM World Trade Corporation, and the IBM Research Mathematical Sciences Department},
  year      = {1972},
  editor    = {Raymond E. Miller and James W. Thatcher and Jean D. Bohlinger},
  pages     = {85--103},
  publisher = {Springer US},
  doi       = {10.1007/978-1-4684-2001-2_9},
  address   = {Boston, MA}
}

@book{vrp3,
  author = {Toth, Paolo and Vigo, Daniele},
  editor = {Daniele Vigo and Paolo Toth},
  title = {Vehicle Routing},
  publisher = {Society for Industrial and Applied Mathematics},
  year = {2014},
  doi = {10.1137/1.9781611973594},
  address = {Philadelphia, PA},
  edition   = {2},
}

@article{ctbrp,
  author    = {Philine Schiewe and Moritz Stinzend{\"o}rfer},
  title     = {Optimizing combined tours: {The} truck-and-cargo-bike case},
  journal   = {OR Spectrum},
  volume    = {46},
  number    = {2},
  pages     = {545--587},
  year      = {2024},
  doi       = {10.1007/s00291-024-00754-2}
}

@article{anderluh,
title = {Multi-objective optimization of a two-echelon vehicle routing problem with vehicle synchronization and ‘grey zone’ customers arising in urban logistics},
journal = {European Journal of Operational Research},
volume = {289},
number = {3},
pages = {940-958},
year = {2021},
issn = {0377-2217},
doi = {10.1016/j.ejor.2019.07.049},
author = {Alexandra Anderluh and Pamela C. Nolz and Vera C. Hemmelmayr and Teodor Gabriel Crainic}
}

@manual{gurobi,
  title        = {Gurobi Optimizer Reference Manual},
  author       = {{Gurobi Optimization, LLC}},
  year         = {2024},
  url          = {https://docs.gurobi.com/current}
}

@misc{Eurostat2016,
  author       = {{Eurostat}},
  title        = {Statistics explained: {Urban Europe}},
  year         = {2024},
  publisher    = {Publications Office of the European Union},
  address      = {Luxembourg},
  note         = {ISSN 2443-8219. Available at
                  \url{https://ec.europa.eu/eurostat/statistics-explained/index.php?title=Urban_Europe}}
}

@article{Crainic2016,
  author    = {Teodor~Gabriel Crainic and
               Nicoletta Ricciardi and
               Giovanni Storchi},
  title     = {Advanced freight transportation systems for congested urban areas},
  journal   = {Transportation Research Part C: Emerging Technologies},
  volume    = {12},
  number    = {2},
  pages     = {119--137},
  year      = {2004},
  doi       = {10.1016/j.trc.2004.07.002},
}

@article{Gendreau1999,
  author    = {Michel Gendreau and
               Gilbert Laporte and
               Christophe Musaraganyi and
               {\'E}ric~D. Taillard},
  title     = {A tabu search heuristic for the heterogeneous fleet vehicle
               routing problem},
  journal   = {Computers \& Operations Research},
  volume    = {26},
  number    = {12},
  pages     = {1153--1173},
  year      = {1999},
  doi       = {10.1016/S0305-0548(98)00100-2}
}

@article{Solomon1987,
  author    = {Marius~M. Solomon},
  title     = {Algorithms for the vehicle routing and scheduling problems
               with time window constraints},
  journal   = {Operations Research},
  volume    = {35},
  number    = {2},
  pages     = {254--265},
  year      = {1987},
  doi       = {10.1287/opre.35.2.254}
}

@article{Desrochers1992,
  author    = {Martin Desrochers and Jacques Desrosiers and Marius Solomon},
  title     = {A new optimization algorithm for the vehicle routing problem
               with time windows},
  journal   = {Operations Research},
  volume    = {40},
  number    = {2},
  pages     = {342--354},
  year      = {1992},
  doi       = {10.1287/opre.40.2.342}
}

@article{Ropke2006,
  author    = {Stefan Ropke and David Pisinger},
  title     = {An adaptive large neighborhood search heuristic for the
               pickup and delivery problem with time windows},
  journal   = {Transportation Science},
  volume    = {40},
  number    = {4},
  pages     = {455--472},
  year      = {2006},
  doi       = {10.1287/trsc.1050.0135}
}

@article{Liu1999,
  author    = {Fu-Hung Liu and Sheng-Yuan Shen},
  title     = {The fleet size and mix vehicle routing problem with time
               windows},
  journal   = {Journal of the Operational Research Society},
  volume    = {50},
  number    = {7},
  pages     = {721--732},
  year      = {1999},
  doi       = {10.1057/palgrave.jors.2600763}
}

@article{Salhi2014,
  author    = {Sa{\"i}d Salhi and Arif Imran and Niaz~A. Wassan},
  title     = {The multi-depot vehicle routing problem with heterogeneous vehicle fleet: {Formulation and a variable neighborhood search implementation}},
  journal   = {Computers \& Operations Research},
  volume    = {52},
  pages     = {315--325},
  year      = {2014},
  doi       = {10.1016/j.cor.2013.05.011}
}

@article{Kocc2016,
  author    = {{\c{C}}a{\u{g}}r{\i} Ko{\c{c}} and Tolga Bekta{\c{s}} and
               Ola Jabali and Gilbert Laporte},
  title     = {The fleet size and mix pollution-routing problem},
  journal   = {Transportation Research Part B: Methodological},
  volume    = {70},
  pages     = {239--254},
  year      = {2014},
  doi       = {10.1016/j.trb.2014.09.008}
}

@article{Perboli2011,
  author    = {Guido Perboli and Roberto Tadei and Daniele Vigo},
  title     = {The Two-Echelon Capacitated Vehicle Routing Problem: Models and Math-Based Heuristics},
  journal   = {Transportation Science},
  volume    = {45},
  number    = {3},
  pages     = {364--380},
  year      = {2011},
  doi       = {10.1287/trsc.1110.0368}
}

@article{MurrayC2015,
  author    = {Chase~C. Murray and Amanda~G. Chu},
  title     = {The flying sidekick traveling salesman problem: {O}ptimization
               of drone-assisted parcel delivery},
  journal   = {Transportation Research Part C: Emerging Technologies},
  volume    = {54},
  pages     = {86--109},
  year      = {2015},
  doi       = {10.1016/j.trc.2015.03.005}
}

@article{Agatz2018,
  author    = {Niels Agatz and Paul Bouman and Marie Schmidt},
  title     = {Optimization approaches for the traveling salesman problem
               with drone},
  journal   = {Transportation Science},
  volume    = {52},
  number    = {4},
  pages     = {965--981},
  year      = {2018},
  doi       = {10.1287/trsc.2017.0791}
}

@article{Drexl2012,
  author    = {Michael Drexl},
  title     = {Synchronization in vehicle routing---{A} survey of {VRPs}
               with multiple synchronization constraints},
  journal   = {Transportation Science},
  volume    = {46},
  number    = {3},
  pages     = {297--316},
  year      = {2012},
  doi       = {10.1287/trsc.1110.0400}
}

@misc{openstreetmap,
  author       = {{OpenStreetMap contributors}},
  title        = {{OpenStreetMap}},
  year         = {2025},
  howpublished = {\url{https://www.openstreetmap.org}},
  note         = {Accessed: 2026}
}

\end{document}